\documentclass [11pt,a4paper]{article}

\usepackage[latin1]{inputenc}

\usepackage[dvips]{graphicx}
\usepackage{amsmath}
\usepackage{amsfonts}
\usepackage{amssymb}
\usepackage{graphics}
\usepackage{rotating}
\usepackage{pifont}
\usepackage{wasysym}
\usepackage{theorem}
\usepackage{array}
\usepackage[arrow, curve, matrix]{xy}
\usepackage[top= 35mm, bottom=30mm, left=35mm, right=35mm]{geometry}
\usepackage{color}
\usepackage{url}

\newtheorem{Satz}{Satz}
\newtheorem{Thm}[Satz]{Theorem}
\newtheorem{Lemma}[Satz]{Lemma}
\newtheorem{Cor}[Satz]{Corollary}

\newtheorem{Sum}[Satz]{Summary}

{\theorembodyfont{\rmfamily} \newtheorem{Def}[Satz]{Definition}
{\theorembodyfont{\rmfamily} 

}
\newtheorem{RD}[Satz]{Remark \& Definition}}

\begin{document}

\title{Rational cohomology of $\overline{R}_2$ (and $\overline{S}_2$)}
\author{Sebastian Krug}
\date{}

\maketitle
\pagestyle{headings}

\newcommand{\mc}{\mathcal}
\newcommand{\GO}{\mathcal{O}}
\newcommand{\wh}{\widehat}
\newcommand{\m}{\mbox}
\newcommand{\mbb}{\mathbb}

\newcommand{\doi}{d_0'}
\newcommand{\doii}{d_0''}
\newcommand{\dor}{d_0^r}
\newcommand{\di}{d_1}
\newcommand{\dii}{d_{1:1}}
\newcommand{\Doi}{D_0'}
\newcommand{\Doii}{D_0''}
\newcommand{\Dor}{D_0^r}
\newcommand{\Di}{D_1}
\newcommand{\Dii}{D_{1:1}}

\newcommand{\aop}{\alpha_0^+}
\newcommand{\bop}{\beta_0^+}
\newcommand{\aip}{\alpha_1^+}
\newcommand{\bip}{\beta_1^+}

\newcommand{\aom}{\alpha_0^-}
\newcommand{\bom}{\beta_0^-}
\newcommand{\aim}{\alpha_1^-}

\newcommand{\Aop}{A_0^+}
\newcommand{\Bop}{B_0^+}
\newcommand{\Aip}{A_1^+}
\newcommand{\Bip}{B_1^+}

\newcommand{\Aom}{A_0^-}
\newcommand{\Bom}{B_0^-}
\newcommand{\Aim}{A_1^-}

\newcommand{\bS}{\overline{S}}
\newcommand{\bR}{\overline{R}}
\newcommand{\bM}{\overline{M}}
\newcommand{\bQ}{\overline{Q}}

\newcommand{\toi}{\tau_0'}
\newcommand{\toii}{\tau_0''}
\newcommand{\tor}{\tau_0^r}
\newcommand{\ti}{\tau_1}
\newcommand{\tii}{\tau_{1:1}}

\newcommand{\roa}{\rho_0^{\alpha}}
\newcommand{\rob}{\rho_0^{\beta}}
\newcommand{\ria}{\rho_1^{\alpha}}
\newcommand{\rib}{\rho_1^{\beta}}

\newcommand{\hoa}{\eta_0^{\alpha}}
\newcommand{\hob}{\eta_0^{\beta}}
\newcommand{\hia}{\eta_1^{\alpha}}

\newcommand{\blue}{\textcolor{blue}}

\sloppy In this text we compute the rational cohomology ring of
$\overline{R}_2$, the moduli space of Prym curves of genus 2, which
is, as we also show, isomorphic to the rational Chow ring of this
space. G. Bini and C. Fontanari did the same for $\overline{S}_2$,
the moduli space of spin curves of genus 2, in \cite{MR2526309}. In
computing the cohomology of $\bR_2$ we follow their approach in
large parts, but also have to apply an idea of Orsola Tommasi,
explained below, to compute additional relations in the rational
cohomology ring. We also correct some errors made in
\cite{MR2526309}. Among other things, some of the relations in the
cohmology rings computed there are not correct, and we apply the
idea just mentioned also to $\bS_2$ in order to replace those
relations. We treat the moduli spaces $\bR_2$ of genus $2$ Prym
curves, $\bS_2^+$ of even genus $2$ spin curves, and $\bS_2^-$ of
odd spin curves parallely. One fact we make intensive use of in our
calculations is that all three moduli spaces are isomorphic to
different of what we call moduli spaces of genus $0$ curves with $6$
\emph{partitioned} marked points (c.f. Lemma \ref{2}). Finite surjective morphisms from $\bM_{0,6}$, 
the moduli space of stable genus $0$ curves with $6$ \emph{ordered} marked points, to all of
the three moduli spaces examined in this text exist, and where introduced
in \cite{MR2526309}. They factor in a natural way through the
mentioned isomorphisms from moduli spaces of curves with partitioned marked points. 
We show more generally that the normalization of the locus of hyperelliptic curves in every 
$\bR_g$ and $\bS_g$ is isomorphic to a disjoint union of several moduli spaces of stable genus $0$ curves
with $2g+2$ partitioned marked points.

I would like to thank Orsola Tommasi, who had the idea of using the
morphisms from $\bM_{0,6}$ just mentioned, together with the fact
that they factor through the moduli spaces of curves with partitioned
points, to compute relations in the cohomology rings of the moduli
spaces of our interest.

\newpage
\tableofcontents

\section{Preliminaries} \label{ss1}

In this section we give basic definitions and results, needed in our
text, and fix notation.

Some notation first:

\begin{enumerate}
 \item A curve means a projective one dimensional variety (not necessarily smooth or irreducible, but necessarily reduced).
 \item By the genus of a curve we will always mean the arithmetic genus.
 \item For any ring $B$ and any group $G$ acting on $B$ we denote by $B^G$ the subring of invariants under the action of $G$.
\end{enumerate}

\subsection{Spin- and Prym curves and their moduli spaces.}

\begin{Def}
(i) A \emph{stable curve} $C$ (possibly with marked points) is a
connected curve having only nodes as singularities and having a
finite group of automorphisms (respecting the marked points, if there
are any). Having a finite automorphism group is equivalent to the
following condition: When we consider as ``special points'' on a
irreducible component of $C$ the marked points as well as the points
in which the component meets the rest of $C$, then every component
of genus $0$ must carry at least three special points, and every
component of genus $1$ must carry at least one special point.

(ii) A \emph{semistable curve} $X$ is a connected stable curve
having only nodes as singularities, and such that every connected
component of genus 1 carries at least one special point, and every
component of genus $0$ carries at least two special points.

(iii) A component of genus $0$ of a semistable curve $X$ meeting the rest of $X$ in exactly two points and carrying no marked points is called an \emph{exceptional component} of $X$.

(iv) The \emph{non-exceptional subcurve} $\tilde X$ of a semistable
curve $X$ is the closure of the complement of all exceptional
components of $X$.

(v) A semistable curve $X$ is called \emph{quasistable}, if no two of its exceptional components intersect each other.

(vi) The \emph{stable model} of a quasistable curve is the (unique) stable curve $C$ obtained by contracting every exceptional component of $X$ to a point. The blow  down map
$\beta: X \rightarrow C$ is also called the stable model of $X$.

\end{Def}

\begin{Def} \label{bf}
(i) A \emph{spin curve} resp. \emph{Prym curve} of genus $g$ is a
triple $(X; \mc{L};b)$, where $X$ is a quasistable curve with stable
model $\beta: X \rightarrow C$, $\mc{L}$ is a line bundle on $X$.
For a spin curve, $b$ is a homomorphism $b: \mc{L}^{\otimes 2}
\rightarrow \omega_X$ such that the restriction of $\mc{L}$ to any
exceptional component $E$ is isomorphic to $\mc{O}_E(1)$ and the
restriction of $b$ to the non-exceptional subcurve $\tilde{X}$
induces an isomorphism $\mc{L}_{|\tilde{X}}^{\otimes 2} \rightarrow
\omega_{\tilde X}$. For a Prym curve replace $\omega_X$ by
$\mc{O}_X$ and $\omega_{\tilde X}$ by $\mc{O}_{\tilde X}$ in the
above definition, and additionally forbid the case $\mc{L} \cong
\GO_X$. The curve $X$ is called the \emph{support} of the spin-
resp. Prym curve, the pair $(\mc{L};b)$ a spin- resp. Prym structure
on $X$. A spin- resp. Prym curve is called smooth if $X$ is smooth.

(ii) We use the definition of isomorphisms of spin curves resp. Prym
curves as for example given in \cite{MR1082361} resp.
\cite{MR2639318}. Thus isomorphisms of spin- resp. Prym curves are
for us isomorphisms of the underlying quasistable curve $X$
compatible with the extra structure,  and do not include a morphism
of the extra structure, as for example in \cite{MR1135424} and
\cite{MR2551759}. I.e. a isomorphism $\varphi: (X; \mc{L}; b)
\rightarrow (X'; \mc{L}'; b')$ of spin- resp. Prym curves is a
isomorphism $\varphi: X \rightarrow X'$ such that there is an isomorphism $\varphi^* \, \mc{L}'
\cong \mc{L}$ which is compatible with $b$. This choice of definition influences the number of
automorphisms of the objects.

(iii) For a given quasistable curve $X$ we call every line bundle
(i.e. invertible sheaf) $\mc{L}$ that fits into the definition of a
spin curve or Prym curve with support $X$ a spin sheaf resp. a Prym
sheaf of $X$. \textbf{We also call the trivial sheaf a Prym sheaf,
and speak of nontrivial Prym sheaves if we want to exclude it.}

(iv) Let $(X;\mc{L};b)$, $(X';\mc{L}';b')$ be two spin- or two Prym
curves, Let $C$, $C'$ be the stable models of $X$ resp. $X'$, let
$N$, $N'$ be the sets of nodes of $C$ resp. $C'$, to which
exceptional components are contracted (``exceptional nodes''). Then
there is a surjective homomorphism of isomorphism groups
\[\psi': Isom(X,X') \rightarrow Isom((C;N),(C';N'))\]
which can of course be restricted to a group homomorphisms
\[\psi: Isom( (X;\mc{L};b),(X';\mc{L}';b')) \rightarrow Isom((C;N),(C';N'))\]
The isomorphisms lying in the kernel of $\psi$ are called
\emph{inessential isomorphisms}. In case of
$(X;\mc{L};b)=(X';\mc{L}';b')$ we speak of \emph{inessential
automorphisms}.
\end{Def}

For every $g \ge 2$ there exist coarse moduli spaces $\bS_g$ and
$\bR_g$ for spin curves resp. Prym curves of genus $g$. They are
projective algebraic varieties of dimension $3g-3$ and have only
finite quotient singularities. The open subsets parametrizing smooth
spin- resp. Prym curves are denoted by $S_g$ and $R_g$. $\bS_g$
consists of two connected components $\bS_g^+$ and $\bS_g^-$
prametricing even resp. odd spin curves.

\begin{Def}
We denote by $\pi_R: \bR_2 \longrightarrow \overline{M}_2$, $\pi_+:
\bS_2^+ \longrightarrow \overline{M}_2$ and $\pi_-: \bS_2^-
\longrightarrow \overline{M}_2$ the ``forgetful morphisms'',
which corresponds to discarding the additional Prym or spin structure,
and passing from $X$ to its stable model $C$.
\end{Def}



\textbf{Notation for other moduli spaces used in this text:}

$\bM_{g,n}$, $\bS_{g,n}$, $\bS_{g,n}^+$, and so on, denote the moduli spaces of
genus $g$ stable curves, spin curves, even spin curves, and so on, together with $n$ ordered
marked points on the underlying curve.

$\bS_{g,n}^{(r_1,...,r_n)}$ resp. $\bR_{g,n}^{(r_1,...,r_n)}$, for
$r_1,...,r_n \in \mbb{Z}$, are moduli spaces of \emph{twisted} spin-
resp. Prym curves with $n$ ordered marked points. Such twisted spin
resp. Prym curves are defined varying the definition of a spin-
resp. Prym curve as follows: If $(p_1,...,p_n)$ are the marked
points on $X$, then the line bundle $\mc{L}$ on $X$ is a square root
of $\omega_X(r_1 p_1+....r_n p_n)$ resp. $\mc{O}_X(r_1 p_1+....+r_n
p_n)$, instead of $\omega_X$ resp. $\mc{O}_X$.

\subsection{Cohomology and rational Chow ring for moduli spaces.}

We will work with the rational Chow ring as well as with the
rational cohomology of moduli spaces. We denote them by
$A^{*}_{\mbb{Q}}(...)$ resp. $H^{*}_{\mbb{Q}}(...)$.

We compile some results by J.H.M Steenbrink from \cite{MR0485870}
about the cohomology of what he calls $V$-manifolds, which are what
we would nowadays call the underlying spaces of orbifolds. All
moduli spaces we are concerned with in this text are $V$-manifolds.

\begin{Sum} \label{p1}
Let X be a projective $V$-manifold. Then

(i) The hard Lefschetz theorem holds, i.e.: Let $L \in H^2(X,
\mbb{Z})$ be the cohomology class of an ample divisor on $X$. Then
for all $q \in \mbb{N}$ the map $\omega \mapsto L^q \wedge \omega$
induces an isomorphism between $H^{n-q}(X,\mbb{C})$ and $H^{n+q}(X,
\mbb{C})$. (\cite{MR0485870} Thm. 1.13)

(ii) The canonical Hodge structure of $H^k(X)$, that would be mixed
for an arbitrary singular variety, is pure of weight $k$ for all $k
\ge 0$. (\cite{MR0485870} Cor. 1.5)
\end{Sum}

Part (ii) allows us to speak of the pure Hodge structure on our moduli spaces, and especially to define Hodge numbers.

In \cite{MR717614} D. Mumford introduced the rational Chow ring of
$Q$-varieties and $Q$-stacks. Our moduli spaces are $Q$-stacks (with
smooth global covers) so we can use Mumfords results. We summarize
the ones we will use:

\begin{Sum} \label{p2}
Let $X$ be an algebraic variety that is a $Q$-variety or a
$Q$-stack, with global Cohen-Macaulay cover. then:

(i) There is a ``natural'' way to define an intersection product
$\bullet \, . \, \bullet$ on the rational Chow group of $X$, making
it into the Chow ring $A^*_{\mbb{Q}}(X)$ we are going to use in our
computations. (C.f. \cite{MR717614} section § 3.)

(ii) To a closed codimension $n$ subvariety $Y$ of one of our moduli
spaces, one can assign classes in the rational Chow ring in two
ways. One is the usual of just taking the corresponding cycle class
$[Y]$ in the Chow group $A^n_Q(X)$. The other is to take the
$Q$-class $[Y]_Q$ of $Y$ as defined in \cite{MR717614} § 3. This
corresponds to considering the cycle of $Y$ on the moduli stack.

(iii) For our moduli spaces, between these two classes the relation
$[Y]=n[Y]_Q$ holds, where $n$ is the number of automorphisms of an
object parametrized by a general point of $Y$.

(iv) Intersections of $Q$-classes in the rational Chow ring, can be
computed on smooth sheets $X_{\alpha}$ mapping to dense open parts
of $X$ (c.f. \cite{MR717614} §3). For our moduli spaces, like for
$\bM_g$, these sheets can be taken to be certain moduli spaces
paremetrizing spin- rep. Prym curves together with a kind of level
structure (c.f. \cite{MR717614} § 3). Locally at any point, the
$X_{\alpha}$ are isomorphic to the deformation space of the object
parametrized by this point. On the deformation space each two
subspaces parametrizing curves of two given topological types, meet like subvectorspaces of a vectorspace, since
local coordinates can be choosen such that one coordinate each coresponds to smoothing of one of the nodes of the curve (C.f. \cite{MR1082361} § 5).
Therefore if one has two cycles $Y$, $Z$ parametrizing
generically curves of a given topological type (as will usually be
the case for the cycles appearing in this text), and their
intersection $Y \cap Z = S$ is proper, then one can treat the
intersection also as transversally in computing the intersection of
the $Q$-classes. Thus in this case $[Y]_Q . [Z]_Q = [S]_Q$. 

(v) A morphism $f: X \rightarrow Y$ of $Q$-stacks (with global
Cohen-Macauley cover), induces a pullback $f^*: A^*_{\mbb{Q}}(Y)
\rightarrow A^*_{\mbb{Q}}(X)$ that is a ring homomorphism. If $W$ is
a closed subvariety of $Y$ such that $codim \; f^{-1}(W)= codim \;
W$, and if we denote by $S$ the set of components of $f^{-1}(W)$
then:
\[f^*([W]_Q) = \sum_{V_k \in S} i_k \cdot [V_k]_Q\]
where $i_k$ can be calculated as the the ramification index of the map $f_{\alpha}: X_{\alpha} \rightarrow
Y_{\alpha}$ belonging to $f$, in the
locus corresponding to $V_k$ on one of the smooth sheets
$X_{\alpha}$. As mentioned above, in our cases these
sheets locally are local universal deformation spaces. (C.f.
\cite{MR717614} Section §3., especially Prop. 3.8.)

(vi) For the pullback $f^*$ just introduced and the usual pushforward $f_*$ the projection formula
(also called push-pull formula) holds:
\[f_*(a.f^*b)= f_*a.b\]
For every $a \in A^*_{\mbb{Q}}(X)$ and $b \in A^*_{\mbb{Q}}(Y)$.
\end{Sum}

In section \ref{s2} we will show that our moduli spaces are even
global quotients of a manifold by a finite group $G$, so in our
special case Steenbrink's and Mumford's results could be shown more
easily:

\begin{Lemma} \label{p3}
 Let $X$ be a smooth algebraic variety, let $G$ be a finite group acting algebraically on $X$ and let $Y=X/G$ be the quotient. Then

 (i) $H^{*}_{\mbb{Q}}(Y) = \bigl ( H^{*}_{\mbb{Q}}(X) \bigr )^G$ (C.f. \cite{MR0413144} Page 120.)

 (ii) $A^{*}_{\mbb{Q}}(Y) = \bigl ( A^{*}_{\mbb{Q}}(X) \bigr )^G$ (C.f. \cite{MR1644323}, Example 1.7.6.)
\end{Lemma}

\subsection{Further notation and conventions for this article}

\begin{enumerate}
 \item If we denote a cycle class of a moduli space by $1$ we mean by this the $Q$-Class of the whole space.
 \item We sometimes speak of ``the stratification according to topological type (of the underlying stable curves)''
 of the spaces $\bR_2$, $\bS_2^+$ or $\bS_2^-$. What is meant
by this is explained in the appendix.
 \item We call ``closed strata'' of these stratifications, the closures of all their strata, not only the strata
 that are already closed (i.e. points).
 \item If there appears a cycle class in our computations that is not written as a product of boundary classes,
 it is usually the class of one of the closed strata just mentioned. (For example $[C^+]_Q$, $[X^-]$, $[E^{\prime, \prime}]$.)
 The (closed) strata are described in the appendix, and they will be used in the main body of the article without defining them there.
 \item If $O$ is an object of the kind parametrized by a moduli space $M$, then we denote the point in $M$ prametrizing $O$ as $[O]$.
 For example if $(X;\mc{L}; b)$ is a Prym curve of genus $g$, then $[(X;\mc{L};b)]$ is the corresponding point in $\bR_g$.
 \item Usually instead of $a.b$ we write $ab$ for the intersection of cycle classes $a,b$ in the Chow ring.
\end{enumerate}

\subsection{Some lemmata for extending morphisms} \label{ss4}

We call a morphism of complex analytic spaces \emph{finite} if it is
proper and has finite fibers. The following lemmata can be proven
quite easily using basic theorems form complex analysis and
commutative algebra.

\begin{Lemma} \label{t5}
 Let $X$, $Y$ be complex analytic spaces, $X$ normal, and $U$ a dense open subset of $X$. If $f: U \rightarrow Y$ is a holomorphic map, and $\tilde f: X \rightarrow Y$ is a continous map
extending $f$, then $\tilde f$ is holomorphic.
\end{Lemma}

\begin{Lemma}\label{t3}
(i) Let $X$, $S$ and $M$ be complex analytic spaces, $X$ normal, $U \subset X$ an open subset.
Let $\pi: S \longrightarrow M$ be a finite holomorphic map,
and let $g:X \longrightarrow M$ and $f:U \longrightarrow S$ be holomorphic maps, such that the following diagram commutes:
\[ \begin{xy} \xymatrix{X \ar@^{<-)}[r]\ar[rrd]_g &  U \ar[r]^f & S \ar[d]^{\pi} \\
                          & & M } \end{xy}\]
Then $f$ extends to a holomorphic map $\tilde{f}: X \longrightarrow S$, compatible with the diagram.

(ii) If furthermore $g$ is finite, then $\tilde{f}$ is finite too.
\end{Lemma}

\begin{Lemma} \label{t4}
Let $X$, $Y$ be algebraic varieties, $Y$ normal. Let $f:X
\rightarrow Y$ be a finite morphism of degree $1$, then $f$ is an
isomorphism.
\end{Lemma}

\subsection{The boundary components of $\overline{R}_2$} \label{ss2}

In the following section we quote results from \cite{MR2639318} we are going to use.

We call the irreducible components of $\bR_2 \smallsetminus R_2$ the
\emph{boundary components} of $\bR_2$. There are exactly $5$ such
components. They have codimension $1$, so they are divisors of
$\bR_2$. The boundary divisors of $\bR_2$ lie above the two boundary
divisors $\Delta_0$ and $\Delta_1$ of $\bM_2$, with respect to the
forgetful map $\pi$. We describe the boundary divisors, by
explaining which kind of Prym curves $(C; \mc{L}; b)$ their general
points parametrize.
\begin{enumerate}
 \item $\Di$: Here $C$ has two irreducible components (of genus 1), meeting in one node, such that restricting the Prym sheaf $\mc{L}$
to one of the components yields the trivial sheaf.
 \item $\Dii$: Here $C$ has two irreducible components meeting in one node, and restricting $\mc{L}$ to either component yields a nontrivial Prym sheaf.
 \item $\Doi$: Here $C$ has one node, the normalization $\tilde{C}$ of $C$ is connected and the pullback of $\mc{L}$ to the normalization is a nontrivial Prym sheaf of $\tilde{C}$.
 \item $\Doii$: Here $C$ has one node, the normalization $\tilde{C}$ of $C$ is connected and the pullback of $\mc{L}$ to the normalization is the trivial sheaf $\GO_{\tilde{C}}$
 \item $\Dor$: Here $C$ consists of two irreducible components, one is a smooth genus 1 curve $D$, the other an exceptional component $E$,
i.e. a smooth genus 0 curve meeting $D$ in two points. Restricting $\mc{L}$ to $D$ yields a Prym sheaf on $D$.
If $\tilde{D}$ and $\tilde{E}$ are the two connected components of the normalization $\tilde C$ of $C$, and if $p,q$ are the two points on $\tilde{D}$
lying over the points of $C$ in which $D$ and $E$ meet, then the Pullback of $\mc{L}$ to $\tilde{D}$ is a square root of $\GO_{\tilde D} (-q-p)$.
\end{enumerate}

To the boundary components we assign elements of
$A_{2,{\mbb{Q}}}(\bR_2)$ by taking $Q$-classes:
\[\di := [\Di]_{Q}, \quad \dii:=[\Dii]_{Q}, \quad \doi:=[\Doi]_{Q}, \quad \doii:=[\Doii]_{Q}, \quad  \dor:=[\Dor]_{Q}\]
we often call these the \emph{boundary classes} of $\bR_2$. Equivalently one defines the boundary classes $\delta_0$ and $\delta_1$ of $\bM_2$.

The forgetful map $\pi_R: \bR_2 \longrightarrow \bM_2$, is ramified in codimension $1$ only at $\Dor$ (therefore the r).
The boundary classes of $\bM_2$ pull back to $\bR_2$ as follows:

\[ \pi^{\ast} (\delta_0) = \doi +\doii +2 \dor \qquad \text{and} \qquad \pi^{\ast} (\delta_1) = \di +\dii \]

The boundary components $\Aop, \Bop, \Aip, \Bip$ of $\bS_2^+$ and
$\Aom, \Bom, \Aim$ of $\bS_2^-$ are described in \cite{MR2526309}.
Again we define corresponding classes:

\[\aop := [\Aop]_{Q}, \quad \bop:=[\Bop]_{Q}, \quad \aip:=[\Aip]_{Q}, \quad \bip:=[\Bip]_{Q},\]
\[\aom := [\Aom]_{Q}, \quad \bom:=[\Bom]_{Q}, \quad \aim:=[\Aim]_{Q} \]

The pullbacks of $\delta_0$ and $\delta_1$ to these spaces are:

\[ \pi_+^* (\delta_0)= \aop + 2 \bop, \quad  \pi_+^* (\delta_1)= 2 \aip +2 \bip,\]
\[ \qquad \pi_-^* (\delta_0)= \aom + 2 \bom, \quad  \pi_-^* (\delta_1)= 2 \aim \]

\section{Moduli spaces of stable hyperelliptic spin- and Prym curves} \label{s2}

\begin{Def}
By $HM_g$, $HS_g^+$, $HS_g^-$ and $HR_g$ we denote the loci of
hyperelliptic curves in $M_g$, $S_g^+$, ... The closures of these
loci in $\bM_g$, $\bS_g^+$, $\bS_g^-$ resp. $\bR_g$ we denote by
$\overline{HM}_g$, $\overline{HS}_g^+$, $\overline{HS}_g^-$, resp.
$\overline{HR}_g$. We call those compact spaces moduli spaces of
stable hyperelliptic curves resp. stable hyperelliptic spin/Prym
curves.
\end{Def}

In this section we show that the normalizations of these compact
moduli spaces, are isomorphic to disjoint unions of several of what
we call moduli spaces of stable genus $0$ curves with
\emph{partitioned} marked points (c.f. Definition \ref{d1}). These
moduli spaces can be described as quotients by finite groups acting
on moduli spaces $\bM_{0,2g+2}$ of stable genus $0$ curves with
$2g+2$ \emph{ordered} marked points. The cohomology rings of the
latter moduli spaces are known by work of S. Keel
(\cite{MR1034665}).

To construct the isomorphisms we will use the fact, that for every
set of $2g+2$ distinct points in $\mbb{P}^1$ there is a (unique up
to isomorphism) $2:1$ cover $h: C \longrightarrow \mbb{P}^1$
ramified exactly over the given points, and $C$ is a genus $g$
smooth hyperelliptic curve, and the fact that every hyperelliptic
curve can be obtained in this way. The spin- resp. Prym sheaves on
$C$ can be recovered as the invertible sheafs corresponding to
certain divisors that are linear combinations of the ramification
points. Using admissible $2:1$ covers of stable genus $0$ curves
with $2g+2$ marked points, one can extend this correspondence to the
asserted isomorphisms.

Probably everything proven in this section is somehow known.

\subsection{Admissible (double) covers}

\begin{Def} \label{d1}
By $\bM_{g,(n_1,...,n_l)}$ we denote the coarse moduli space of the
following moduli problem: Stable curves of genus $g$ with
$n_1+n_2+...+n_l$ unordered marked points that are divided into $l$
disjoint sets $A_1,...,A_l$ such that $\# A_i = n_i$ for all $i \in
\{1,...,l\}$. Formally the objects are pairs $(X; (A_1,...,A_l))$,
where $X$ is a genus $g$ curve such that $X$ becomes stable if one
marks the points on it contained in the $A_i$. The moduli space one
gets by changing the objects of the moduli problem to $(X;
\{A_1,...,A_l\})$, i.e. by having a set of sets instead of a tuple
of sets in the defining data, we denote by $\bM_{g, [ n_1,...,n_l
]}$.

We often call moduli spaces of the latter type, moduli spaces of curves of genus $g$ with \emph{partitioned} marked points.

By $M_{g,(n_1,...,n_l)}$ resp. $M_{g, [ n_1,...,n_l ]}$ we denote the corresponding moduli spaces of smooth curves.
\end{Def}

\begin{RD} \label{rd1}
(i) For $n:=n_1+n_2+...+n_l$ one can construct the moduli space
$\bM_{g,(n_1,...,n_l)}$ as the quotient of $\bM_{g,n}$. Divide the
set of marked points $\{1,...,n\}$ into disjoint subsets
$A_1',...,A_l'$ of the appropriate size $\# A_i'= n_i$. Then take
the quotient of $\bM_{g,n}$ induced by the action of $S_{n_1} \times
... \times S_{n_l}$ permuting the indices inside the sets
$A_1',...,A_n'$. $\bM_{g,[n_1,...,n_l]}$ can be constructed as the
quotient of $\bM_{g,(n_1,...,n_l)}$ by the action permuting the
indices of those of the sets $A_1,...,A_n$ having the same
cardinality.

(ii) For genus $0$ we fix some of the quotient morphisms, we are going to use later. Let
\[\pi_{(n_1,...,n_l)}: \bM_{0,n} \rightarrow \bM_{0,(n_1,...,n_l)}
\]
be the quotient morphism corresponding to the choice $A_1':=
\{1,...,n_1\}$, $A_2'=\{n_1+1,n_1+2,...,n_1+n_2\}$, and so on.

Let
\[\pi_{[n_1,...,n_l]}: \bM_{0,n} \rightarrow \bM_{0,[n_1,...,n_l]}
\]
be the composition of $\pi_{(n_1,...,n_l)}$ with the quotient
morphism $\bM_{0,(n_1,...,n_l)} \rightarrow \bM_{0,[n_1,...,n_l]}$.
\end{RD}

\begin{Def} \label{d2}
 (i) Let $(D;\{p_1,...,p_n\})$ be a stable genus $0$ curve with $n$ unordered marked points.
 For us an admissible $d:1$ cover of $(D; \{p_1,...,p_n\})$ is a regular morphism $f:Y \rightarrow D$ such that $Y$ is
 a connected nodal curve, and:
 \begin{enumerate}
 \item For $D_{ns}$ the nonsingular locus of $D$, $f^{-1}(D_{ns})=Y_{ns}$ and the restriction of $f$ to $f^{-1}(D_{ns})$ is a $d:1$ cover simply ramified over the marked points and unramified everywhere else.
 \item For every node $q$ of $D$, every point in $f^{-1}(q)$ is a node of $Y$ and for every such node $r$ the two branches of $Y$ in $r$ are mapped to the two branches of $D$ near $q$,
 both with the same ramification index in $r$.
 \end{enumerate}
 (ii) An isomorphism between two admissible covers $f: Y \rightarrow D$ and $f':Y' \rightarrow D'$ is an isomorphism $\varphi: Y \rightarrow Y'$ such that there is an
isomorphism $\psi: D \rightarrow D'$ for which $\psi \circ f = f'
\circ \varphi$.
\end{Def}

We compile some facts about admissible covers, especially $2:1$
covers, which we mostly take from \cite{MR1981190}.

\begin{Sum} \label{d3}

(i) There is a coarse moduli space $\overline{H}_{d,g}$ of admissible $d:1$-covers of stable genus $0$ curves with $2(g+d)-2$ marked points.

(ii) The covering space $Y$ of an admissible $2:1$ cover $f:Y
\rightarrow D$ is a semistable curve, all whose irreducible
components are smooth.

(iii) There are isomorphisms $\bM_{0,[2g+2]} \cong
\overline{H}_{2,g} \cong \overline{HM}_g$. They are explicitly
constructed in \cite{MR1981190} by describing how to associate to a
family of stable genus $0$ curves with $2g+2$ unordered marked
points a unique family of admissible $2:1$ covers, and how to
associate to a family of such admissible covers a unique family of
stable hyperelliptic curves. (The latter is just by contracting all
exceptional components of $Y$, i.e. stable reduction.)

(iv) Thus, in particular, the $2:1$ admissible cover $f:Y
\rightarrow D$  of a stable genus $0$ curve with an even number of
unordered marked points is defined uniquely up to isomorphism.
\end{Sum}

\textbf{Caution:} Sometimes we will just mean the $Y$ of $f: Y \rightarrow D$ when we talk about the
admissible $2:1$ cover of a stable genus $0$ curve with marked points.

\subsection{Relation to moduli spaces of stable genus $0$ curves with partitioned marked points.} \label{ss5}

\begin{Lemma} \label{a6}
For $g \ge 2$, let $p_1,...,p_{2g+2}$ be distinct points in $\mbb{P}^1$, and $h: Y \longrightarrow \mbb{P}^1$ the (unique) $2:1$
cover of $\mbb{P}^1$ ramified
exactly over these points. Then $Y$ is a genus $g$ hyperelliptic curve.
For $i=1,...,2g+2$, define $q_i := h^{-1}(p_i)$. Let $Q$ be the set of all
$q_i$ and denote by $P_n$ the set of possible partitions of $Q$ into a set of $n$ elements and a set of $2g+2-n$ elements.
I.e.:
 \[P_n:= \{ \{A,B\} \; | \; A,B \subseteq Q, \, A \uplus B =Q, \, \# A = n, \, \# B = 2g+2-n \}\]
Let $J_R(Y)$, $J_S(Y)$, $J_+(Y)$, $J_-(Y)$ be the sets of
isomorphism classes of nontrivial Prym sheaves, resp. spin sheaves,
resp. even spin sheaves, resp. odd  spin sheaves on Y. (Of course
$J_S (Y) = J_+ (Y) \uplus J_- (Y)$.) Then we have:

For any $\{A,B\} \in P_n$ and $r_1,...,r_n$ the points in A.

(i) For all even $2 \le n \le g+1$:
\begin{enumerate}
 \item $\phi_{R,n}(\{A,B\}):= \GO_Y(r_1+...+ r_{\frac{n}{2}}- r_{\frac{n}{2}+1} -...- r_{n})$ is a nontrivial Prym sheaf of $Y$. Its
isomorphism class is independent of the ordering of the points $r_i$, as well as of the choice of $A$, necessary in the case $n = g+1$. Thus the following map is well defined.
 \item The map $ \phi_{R,n}: P_n \rightarrow J_R(Y)$, $\{A,B\} \mapsto \phi_{R,n}(\{A,B\})$ is injective.
 \item The map $\phi_{R}: \biguplus_{\substack{2 \le n \le g+1, \\ n \; \text{even}}} P_n \rightarrow J_R(Y)$, obtained as union of the maps $\phi_{R,n}$ is a bijection.
\end{enumerate}

(ii)  Analogously for spin structures:

\begin{enumerate}
\item If $g$ is even, then for all $0 \le n \le g+1$, with $n$ odd:

$\phi_{S,n}(\{A,B\}):= \GO_Y((g-2) \cdot q_1 + r_1+r_2+...+ r_{\frac{n+1}{2}} - r_{\frac{n+1}{2}+1} -...- r_{n})$ is a spin sheaf of $Y$. 
\item If $g$ is odd, then for all $0 \le n \le g+1$, with $n$ even:

$\phi_{S,n}(\{A,B\}):= \GO_Y(g \cdot q_1+ r_1 + r_2 +...+ r_{\frac{n}{2}}- r_{\frac{n}{2}+1} -...- r_{n})$ is a spin sheaf of $Y$. 

\item In both cases the isomorphism class of $\phi_{S,n}(\{A,B\})$ is independent of the ordering of the points $r_i$ and $q_i$, as well as of the choice of $A$, necessary in the case $n = g+1$. Thus the map $ \phi_{S,n}: P_n \rightarrow J_S(Y)$, $\{A,B\} \mapsto \phi_{R,n}(\{A,B\})$ is well defined. It is injective, and the map $\phi_{S}: \biguplus_{\substack{1 \le n \le g+1, \\ n \; \text{odd}}} P_n \rightarrow J_S(Y)$, obtained as union of the maps $\phi_{S,n}$, is a bijection.
\end{enumerate}

(iii) For every $g \ge 2$ the bijection $\phi_S$ of course splits into two bijections $\phi_+: (\phi_S)^{-1} J_+(Y) \rightarrow J_+(Y)$ and
$\phi_-: (\phi_S)^{-1} J_-(Y) \rightarrow J_-(Y)$. They can also be written (by describing $(\phi_S)^{-1} J_+(Y)$ and $(\phi_S)^{-1} J_-(Y)$ explicitly) as:
\[\phi_{+}: \biguplus_{\substack{1 \le n \le g+1, \\ n \equiv g+1 \; \text{mod} \; 4}} P_n \rightarrow J_+(Y)\]
and
\[\phi_{-}: \biguplus_{\substack{1 \le n \le g+1, \\ n \equiv g-1 \; \text{mod} \; 4}} P_n \rightarrow J_-(Y)\]
\end{Lemma}

\textbf{Proof:} It is easy to show that, for all $i,j \in \{1,...,2g+2\}$, $2q_i-2q_j \sim 0$. I.e. all $2 q_i$ are equivalent.

Using this, all claims of part (i) follow form what is shown in § 5.2.3. in \cite{Dolga1}. 

All assertions of (ii) follow from the fact that
the canonical sheaf of $Y$ is equivalent to $(2g-2) q_i$ for any $i \in \{1,...,2g+2\}$ and the corresponding assertions of part (i) of the Lemma.

For (iv): From Lemma 5.2.1. in \cite{Dolga1} it follows that $h^0(\phi_{S,n}(\{A,B\}))$ is even if $g-n+1 \equiv 0 \mod 4$ and odd if $g-n+1 \equiv 2 \mod 4$. This proves part (iv) of the Lemma. $\square$

\begin{Lemma} \label{n1}
If by $X^{\sim}$ we denote the normalization of a variety $X$ then: 

 (i) For all $g \ge 2$ there is an isomorphism:
 \[b: \bM_{0,[2g+2]} \overset{\cong}{\longrightarrow} \overline{HM}_g\]

 (ii) For all $g \ge 2$ there is an isomorphism:
 \[ a_R: \biguplus_{\substack{2 \le n \le g+1, \\ n \; \text{even}}} { \bM_{0,[n,2g+2-n]}} \overset{\cong}{\longrightarrow}  (\overline{HR}_g)^{\sim}\]

 (iii) For all $g \ge 2$ there are isomorphisms:

 \[ a_+: \biguplus_{\substack{0 \le n \le g+1, \\ n \equiv g+1 \; \text{mod} \; 4}} { \bM_{0,[n,2g+2-n]}}  \overset{\cong}{\longrightarrow} (\overline{HS}_g^+)^{\sim}\]
 and
 \[ a_-: \biguplus_{\substack{0 \le n \le g+1, \\ n \equiv g-1 \; \text{mod} \; 4}} { \bM_{0,[n,2g+2-n]}} \overset{\cong}{\longrightarrow}  (\overline{HS}_g^-)^{\sim}\]

All the isomorphism above map boundary points to boundary points (after composing the isomorphisms with the normalization map).
\end{Lemma}

\textbf{Proof:} (i) The isomorphism of (i) is constructed in
\cite{MR1981190} (also c.f. Summary \ref{d3}), it maps boundary
points to boundary points, as can easily be checked by looking at
the construction there.

On the interior of the moduli spaces the restricted morphism $b':
M_{0,[2g+2]} \longrightarrow HM_g$ acts in the following way: Let
$(D; \{p_1,...,p_{2g+2}\})$ be a smooth rational curve with $2g+2$
unordered marked points. Let $Y$ be the unique $2:1$ cover of $D$ ramified
over exactly the points $p_i$. The morphism
$b'$ assigns to $[(D; \{p_1,...,p_{2g+2}\}] \in M_{0,2g+2}$ the
point $[Y] \in HM_g$. Every smooth hyperelliptic curve $Y$ of genus
$g$ is a $2:1$ cover of $\mbb{P}^1$ ramified in $2g+2$ Points, thus
$b'$ is surjective. Since two smooth pointed curves $(D;
\{p_1,...,p_{2g+2}\})$ and $(D'; \{p'_1,...,p'_{2g+2}\})$ are
isomorphic if and only if the covers $Y$ and $Y'$ are isomorphic,
$b'$ is of degree $1$. Both $M_{0,[2g+2]}$ and $M_g$ are normal
varieties, thus this implies that $b'$ is an isomorphism. For the
description of the isomorphism $b$, extending $b'$ to the compactified moduli spaces, c.f.
\cite{MR1981190}.

Now we prove (ii): A morphism
\[a_R': \biguplus_{\substack{2 \le n \le g+1, \\ n \; \text{even}}} { M_{0,[n,2g+2-n]}} \longrightarrow HR_g\]
can be defined in the following way: For $2 \le n \le g+1$, $n$
even, $a'_R$ assigns to $[(D; \{A,B\})] \in M_{0,[n,2g+2-n]}$ the
point $[(Y; \phi_{R,n}(\{A,B\}))] \in HR_g$, where $Y$ is defined as
above, and $\phi_{R,n}(\{A,B\})$ is defined as in Lemma \ref{a6}
(i). The morphism $a'_R$ is surjective and $1:1$ by Lemma \ref{a6}
(i). Let
\[ \psi: \biguplus_{\substack{2 \le n \le g+1, \\ n \; \text{even}}} \bM_{0,[n,2g+2-n]} \longrightarrow \bM_{0,[2g+2]} \]
and
\[ \pi: \overline{HR}_g \longrightarrow \overline{HM}_g \]
be the forgetful morphisms. Using the abbreviations $N:= \biguplus_{\substack{2 \le n \le g+1, \\ n \; \text{even}}} M_{0,[n,2g+2-n]}$ and
$\overline N := \biguplus_{\substack{2 \le n \le g+1, \\ n \; \text{even}}} \bM_{0,[n,2g+2-n]}$, we get the commutative diagram

\[ \begin{xy} \xymatrix{\overline N \ar@^{<-)}[r]\ar[rrd]_{b \circ \phi}
&  N \ar[r]^{a_R'} & \overline{HR}_g \ar[d]^{\pi} \\
                          & & \overline{HM}_g } \end{xy}\]

Thus by Lemma \ref{t3} $a_R'$ extends to a finite  surjective Morphism
\[a_R: \biguplus_{\substack{2 \le n \le g+1, \\ n \; \text{even}}} \bM_{0,[n,2g+2-n]} \longrightarrow \overline{HR}_g.\]
It has degree $1$ and must be an isomorphism since both varieties are normal (c.f. Lemma \ref{t4}).

Part (iii) of our Lemma is proven analogously to part (ii), by using part (ii) and (iii) of Lemma \ref{a6} instead of part (i).
The isomorphisms of part (ii) and (iii) of our Lemma map boundary points to boundary points because they are compatible
with the isomorphism $b$ of Part (i), and this one does.
$\square$

\textbf{Remark:} While $\overline{HM}_g$ is a normal variety for all $g \ge 2$ (since it is isomorphic to $\bM_{0,[2g+2]}$), the spaces $\overline{HS}_g^+$, $\overline{HS}_g^-$ and
$\overline{HR}_g$ in gereral are not. Take for example in $\bS_3^-$ a point coresponding to a spin curve $(X;\mc{L};b)$ with $X$ consisting of two disjoint smooth genus $1$ curves and two exceptional components, such that each exceptional component meets each genus $1$ component in exactly one point. Now let $D'$ be the local universal deformation space of $(X;\mc{L};b)$, and let $D$ be the local universal deformation space of $C$, the stable model of $X$. Then the forgetful map $\varphi: D' \rightarrow D$ is $4:1$ and simply ramified over each of the two subspaces of $D$ coresponding to the two nodes of $C$, which are blown up in $X$. One can define local coodinates $x,y,z_1,...,z_4$ around the special point of $D$, such that the two subspaces just mentioned are the spaces $x=0$ and $y=0$. Then, for suitably choosen local coordinates on $D$,  $\varphi$ is described by $x' \mapsto (x')^2$, $y' \mapsto (y')^2$ and $z_i' \mapsto z_i'$ for $i=1,...,n$. The hyperelliptic involution on $C$ swaps the two nodes, thus the hyperelliptic locus in $D$ is invariant under swaping the coordinates $x$ and $y$. Since the hyperelliptic locus is also normal, it follows that it can localy be described by $x=y$ (for a possible choice of coordinates). Thus the hyperelliptic locus in $D'$ is described by $(x')^2=(y')^2$, thereby having a singularity of codimension $1$. As one can check, this singularity is retained when quotienting $D'$ by the action of the automorphism group of $(X;\mc{L};b)$, hence the hyperelliptic locus $\overline{HS}_3^-$ in $\bS_3^-$ is not normal.     

\subsection{Some properties of $\bM_{0,n}$}

The moduli spaces $\bM_{0,n}$ ($n \ge 3$) of stable genus $0$ curves
with ordered marked points where examined by S. Keel in
\cite{MR1034665}. Among other things he computed their cohomology
ring (and, what is the same for these spaces, the Chow ring) for all
$n \ge 3$. We summarize some facts about these spaces we are going
to use from \cite{MR1034665}.

\begin{Sum} \label{n2} \textbf{\emph{(S. Keel)}}

For all $n \ge 3$:

(i) $\bM_{0,n}$ is a smooth rational projective variety of dimension $n-3$.

(ii) For every $S \subsetneq \{1,...,n\}$ such that $\# S \ge 2$ and
$\# (\{1,...,n\}\smallsetminus S) \ge 2$, there is a boundary
divisor $D^S$ of $\bM_{0,n}$, a general point of which corresponds
to a rational curve with two smooth irreducible components meeting
in one node, such that the marked point with indices in $S$ lie on
one of the components, and the marked points with indices in $S^c:=
\{1,...,n\}\smallsetminus S$ lie on the other component. (Of course
$D^S$ and $D^{S^c}$ are the same divisor.)

The boundary $\bM_{0,n} \smallsetminus M_{0,n}$ of $\bM_{0,n}$ is exactly the union of the divisors just described.

(iii) The cohomology ring of $\bM_{0,n}$ is generated by the
boundary components, and is isomorphic to the chow ring by the cycle
map.

(iv) More specific:
\[H^*(\bM_{0,n}) \cong A^*(\bM_{0,n}) = \frac{\mbb{Z}[ \{ D^S |  S \subsetneq \{1,...,n\}, \; \# S \ge 2, \;\# S^c \ge 2 \} ]}{
\{\text{the following relations}\}}\]
The relations in the Chow ring are:

\begin{enumerate}
 \item For all $S \subsetneq \{1,...,n\}$ such that $\# S \ge 2$ and $\# S^c \ge 2$: $D^S = D^{S^c}$
 \item For every for $i,j,k,l \in \{1,...,n\}$:
  \begin{equation} \label{a9}
    \sum_{\substack{S \subsetneq \{1,...,n\},\\ i,j \in S,\\ k,l \notin S}} D^S = \sum_{\substack{S \subsetneq \{1,...,n\},\\ i,k \in S,\\ j,l \notin S}}
    D^S =\sum_{\substack{S \subsetneq \{1,...,n\},\\ i,l \in S,\\ j,k \notin S}} D^S
  \end{equation}
 \item For all $S,T \subsetneq \{1,...,n\}$ such that $\# S, \, \# T, \, \# S^c, \, \# T^c \ge 2$: $D^SD^T=0$ if not one of the following conditions holds:
 \[S \subseteq T, \quad T \subseteq S, \quad S \subseteq T^c, \quad S^c \subseteq T \]
\end{enumerate}

\end{Sum}

\begin{Def}\label{short}
We use the following short notation for the boundary components of
$\bM_{0,n}$: If $\{a_1,...,a_l\}$, $2 \le l \le n-2$, is a subset of
$\{1,...,n\}$ we will denote the corresponding boundary divisor
$D^{\{a_1,...,a_l\}}$ by $[a_1,...,a_l]$.
\end{Def}

\subsection{Conclusions}

\begin{Cor} \label{n3}
For all $g \ge 2$ and every $\overline{X} \in \{ \overline{HM}_g, (\overline{HS}_g^+)^{\sim}, (\overline{HS}_g^-)^{\sim}, (\overline{HR}_g)^{\sim} \}$ we have:

(i) Every connected component of $\overline{X}$ is unirational.

(ii) $A_{\mbb{Q}}^*(\overline{X}) \cong H_{\mbb{Q}}^*
(\overline{X})$, as graded $\mbb{Q}$-algebras, after adjusting the
grading of $A_{\mbb{Q}}^*(\overline{X})$ by a factor $2$. In
particular $H_{\mbb{Q}}^n (\overline{X})=0$ for all odd $n$.

(iii) $Pic_{\mbb{Q}}(\overline{X}) \cong A^1_{\mbb{Q}}(\overline{X})$

(iv) $A^1_{\mbb{Q}}(\overline{X})$ is generated by the boundary divisors of $\overline{X}$. 
(Meaning the preimages of the boundary components of the moduli space on its normalization.) 

(v) $h^{p,0}(\overline{X})=0$ for $p > 0$.
\end{Cor}

\textbf{Proof:} For all claims it suffices to show them for every
connected component of $\overline{X}$. Let $\overline{Y}$ be such a
component, $Y$ its Interior. Then, by Lemma \ref{n1} and the Remark
\ref{rd1}, $\overline{Y} \cong \bM_{0,2g+2}/G$ for some subgroup G
of $S_{2g+2}\times S_2$.

(i):  $\overline Y \cong \bM_{0,2g+2}/G$ is of course covered by
$\bM_{0,2g+2}$, and all spaces $\bM_{0,n}$ are rational (Summary
\ref{n2} (i)).

(ii): By Summary \ref{n2} (iii), $A^*_{\mbb{Q}}(\bM_{0,2g+2}) \cong H^*_{\mbb{Q}}(\bM_{0,2g+2})$. Using Lemma \ref{p3} we get:
\[A^*_{\mbb{Q}}(\overline Y) \cong A^*_{\mbb{Q}}(\bM_{0,2g+2}/G) \cong (A^*_{\mbb{Q}}(\bM_{0,2g+2}))^G \]
\[\cong (H^*_{\mbb{Q}}(\bM_{0,2g+2}))^G \cong H^*_{\mbb{Q}}(\bM_{0,2g+2}/G) \cong H^*_{\mbb{Q}}(\overline Y)\]

(iii): $\overline Y$ is normal, so the Picard group is in a natural way a subgroup of the divisor class group, c.f. \cite{MR0463157}
Remark 6.11.2. and Prop. 6.15. Thus there is an injection
\[Pic_{\mbb{Q}}(\overline Y) \longrightarrow A_{\mbb{Q}}^1(\overline Y) \]
Since $\overline Y \cong \bM_{0,2g+2}/G $ has only finite quotient
singularities, it is $\mbb{Q}$-factorial, i.e. every Weil-divisor is
$\mbb{Q}$-Cartier. Thus the map is also surjective.

(iv): By Summary \ref{n2}, $A^1(\bM_{0,2g+2})=
A_{(2g-1)-1}(\bM_{0,2g+2})$ is generated by the boundary classes,
i.e. the map $A_{(2g-1)-1}(\bM_{0,2g+2} \smallsetminus M_{0,2g+2})
\longrightarrow  A_{(2g-1)-1}(\bM_{0,2g+2})$ is surjective. The exact
sequence
\[ A_{(2g-1)-1}(\bM_{0,2g+2} \smallsetminus M_{0,2g+2}) \longrightarrow  A_{(2g-1)-1}(\bM_{0,2g+2}) \longrightarrow A_{(2g-1)-1}(M_{0,2g+2}) \longrightarrow 0 \]
then yields $A_{\mbb{Q}, (2g-1)-1}(M_{0,2g+2}) =
A_{(2g-1)-1}(M_{0,2g+2})=0$. By Lemma \ref{p3}, then
\[A_{\mbb{Q}, (2g-1)-1}(Y) \cong A_{\mbb{Q},(2g-1)-1}(M_{0,2g+2}/G) \cong (A_{\mbb{Q},(2g-1)-1, }(M_{0,2g+2}))^G=0\]
Again using an exact sequence like the one above we conclude that
$A_{\mbb{Q},(2g-1)-1}(\overline{Y} \smallsetminus Y) \longrightarrow
A_{\mbb{Q}, (2g-1)-1}(\overline{Y})$ is surjective, i.e. that
$A_{\mbb{Q},(2g-1)-1}(\overline{Y}) \cong
A_{\mbb{Q}}^1(\overline{Y})$ is generated by the boundary classes.

(v): According to \cite{MR1034665}, every $\bM_{0,2g+2}$, is rational. Thus $H^{p,0}(\bM_{0,2g+2}) \cong H^{p,0}(\mbb{P}^{n-3})=0$ for all $p > 0$, since all $h^{p,0}$ are birational invariants (c.f. \cite{MR1288523} p. 494).
This implies $H^{p,0}(\overline{Y}) = (H^{p,0}(\bM_{0,2g+2}))^G =0 $. $\square$

\section{Morphisms to $\bS_2$ and $\bR_2$.}

In this section we introduce several finite morphisms from other
moduli spaces to $\bR_2$, $\bS_2^+$ and $\bS_2^-$. They will later
be used to determine relations between cycle classes on our moduli
spaces, by pushing forward known relations, or by using push-pull
formula.

\subsection{The morphisms $a_R$, $a_+$ and $a_-$ in the case of genus
2} \label{sec1}

In the case of genus $2$, all smooth curves are hyperelliptic, hence
$\overline{HM}_2= \bM_2$, $\overline{HR}_2= \bR_2$,
$\overline{HS}^+_2= \bS^+_2$ and $\overline{HS}^-_2= \bS^-_2$. Thus
the conclusions listed in Corollary \ref{n3}, apply to the moduli
spaces we are interested in. Lemma \ref{n1} in this special case
reads

\begin{Lemma} \label{2} \textbf{($\&$ Definition)}

 There are Isomorphisms
 \[ b: \bM_{0,[6]} \overset{\cong}{\longrightarrow} \bM_2 \quad \text{resp.}\]
 \[a_R:\bM_{0,[2,4]} \overset{\cong}{\longrightarrow} \bR_2 \quad \text{resp.} \quad a_+: \bM_{0,[3,3]} \overset{\cong}{ \longrightarrow} \bS_2^{+} \quad \text{resp.} \quad a_-: \bM_{0,[1,5]} \overset{\cong}{\longrightarrow} \bS_2^{-} \]
 These isomorphisms map the boundary of $\bM_{0,6}$ onto the boundary of the images.

 We define surjective finite morphisms, by composing every one of isomorphisms above with the appropriate quotient morphism out of,
 $\pi_{0,[6]})$, $\pi_{0,[2,4]}$, $\pi_{0,[3,3]}$, and $\pi_{0,[1,5]}$ from Definition \ref{rd1}:
 \[ g: \bM_{0,6} \overset{720:1}{\longrightarrow} \bM_2, \]
 \[f_R: \bM_{0,6} \overset{48:1}{\longrightarrow} \bR_2, \quad f_+: \bM_{0,6} \overset{72:1}{\longrightarrow} \bS_2^{+} \quad \text{and} \quad f_-: \bM_{0,6} \overset{120:1}{\longrightarrow} \bS_2^{-}. \]
\end{Lemma}

\textbf{Proof:} Everything except the degrees of the finite
surjective morphisms is just a special case of Lemma \ref{n1}, since all genus 2 curves are hyperelliptic, and since our moduli spaces are normal. The degrees equal those of the forgetful morphisms $\pi_{0,[6]})$,
$\pi_{0,[2,4]}$, $\pi_{0,[3,3]}$, $\pi_{0,[1,5]}$, which can easily
be counted. $\square$

By the Lemma we even know:

\[H^*_{\mbb{Q}}(\bR_2) \cong (H^*_{\mbb{Q}}(\bM_{0,6}))^{S_2 \times S_4}, \quad H^*_{\mbb{Q}}(\bS^+_2) \cong (H^*_{\mbb{Q}}(\bM_{0,6}))^{S_3 \times S_3 \times S_2},\]
\[ H^*_{\mbb{Q}}(\bS^-_2) \cong (H^*_{\mbb{Q}}(\bM_{0,6}))^{S_1 \times S_5} \]

where the group actions are those of Remark \ref{rd1}. As the
cohomology of $\bM_{0,6}$ is known (c.f. Summary \ref{n2}), probably
a computer algebra program could compute these invariant subrings.
In this article we instead compute the rational cohomology of
$\bR_2$ and $\bS_2$ by hand. For our computation we need some more
information about the isomorphism $a_R$, $a_+$ and $a_-$, and the
finite surjective maps $f_R$, $f_+$, and $f_-$ defined from them.

First we examine which boundary components are identified by the
isomorphisms $b$, $a_R$, $a_+$ and $a_-$. On $\bM_{0,[6]}$ there are
two boundary divisors corresponding to the two types of smooth
curves genus $0$ curves with $6$ unordered marked points and one
node, which are:
\begin{enumerate}
 \item Two smooth genus $0$ curves meeting in one node, with $4$ marked points on one curve, $2$ marked points on the other.
 \item Two smooth genus $0$ curves meeting in one node, with $3$ marked points on each curve.
\end{enumerate}

From the theory of admissible covers it is clear that the first
boundary divisor is mapped to $\Delta_0$ and the second to
$\Delta_1$ by $b$. If we describe boundary divisors of $\bM_{0,[6]}$
by diagrams of the objects corresponding to general points, we get
the table

\begin{center}
\begin{tabular}{|m{4cm}|m{4cm}|}
\hline Bound. Div. of $\bM_{0,[6]}$ & is mapped to \\
\hline \hline \includegraphics[width=2cm]{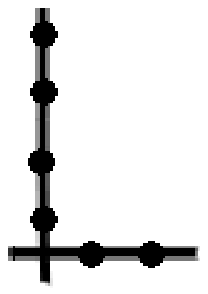} & $\Delta_0$\\
\hline \includegraphics[width=2cm]{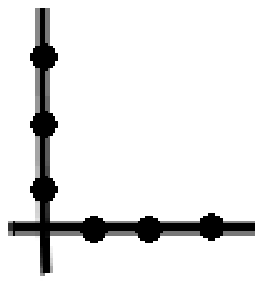} & $\Delta_1$\\
\hline
\end{tabular}
\end{center}

For the isomorphisms $a_R$, $a_+$ and $a_-$ we use that there is the following
commutative diagram for $a_R$, and analogous diagrams for $a_+$ and $a_-$.
\[ \begin{xy}
     \xymatrix{ \bM_{0,[2,4]} \ar[d]^{\psi} \ar[r]^{a_R} & \bR_2 \ar[d]^{\pi_R}
             \\ \bM_{0,[6]} \ar[r]^{b} & \bM_2 }
    \end{xy}  \]
Here $\psi: \bM_{0,[2,4]} \rightarrow \bM_{0,[6]}$ is the forgetful morphism.

Let $C$ be a boundary component of $\bM_{0,[2,4]}$ and $D$ the
boundary component of $\bR_2$ it is mapped to by $a_R$. Take a
general point $y$ in the boundary component of $\bM_2$ underlying
$D$, and let $x:= b^{-1}(x)$ be its preimage in $\bM_{0,[6]}$. Then
the number $M$ of elements in the fiber $\psi^{-1} (x)$ lying in $C$
must be equal to the number $N$ of elements of the fiber $\pi_R^{-1}
(y)$ lying in $D$. Knowing the numbers $N$ and $M$ for all boundary
components, and the behavior of $b$ suffices to see which components
get identified. The number $M$ can be counted if we draw for every
boundary component of $\bM_{0,[2,4]}$ the diagram of the objects
corresponding to a general point. In those diagrams we will always
denote a elements of $A$ by squares and elements of $B$ by dots. (In
the case of $\bM_{0,[3,3]}$ we do not distinguish between $A$ an
$B$, and only say that squares and dots belong to different sets.)

For example

\begin{center} \includegraphics[width=2cm]{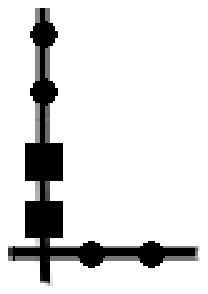} \end{center}

describes an object $(X; \{ A,B\})$ where $X$ consists of two smooth
genus $0$ curves meeting in one node, $A$ is a set of two points on
$X$, $B$ a set of $4$ points on $X$, $A$ and $B$ disjoint, such that
one of the two smooth curves contains two elements of $B$ and all
two elements of $A$, while the other curve contains the remaining
two elements of $B$. One can count that there are $M=\binom{4}{2}=6$
objects of this type lying over the corresponding object in
$\bM_{0,[6]}$ described by the diagram

\begin{center} \includegraphics[width=2cm]{del_0.eps} \end{center}

The numbers $N$ can be determined, using the descriptions of the
boundary divisors in section \ref{ss2}. For example for $\Doi$ the
number is $N =6$. Indeed, if $(X; \mc{L}; b)$ is a general object
parametrized by $\Doi$, then $\mc{L}$ can be obtained by taking a
nontrivial Prym sheaf on the normalization of $X$, and gluing the
sheaves fibers over the points identified by the normalization map,
in one of two possible ways. The Normalization is an elliptic curve,
so there are are $3$ nontrivial Prym sheaves on it, and gluing them
in the two possible ways yields 6 nonisomorphic Prym sheaves on $X$.

After computing $M$ and $N$ for all boundary components one can
conclude that the two boundary components from our examples get
identified. The following tables list the identifications for all
boundary components, together with the numbers $N$ and $M$

\begin{center}
\begin{tabular}{|m{4cm}|m{3cm}|m{3cm}|}
\hline Bound. Div. of $\bM_{0,[2,4]}$ & is mapped to& $M=N$ \\
\hline \hline \includegraphics[width=1.6cm]{doi.eps} & $\Doi$ & $6$ \\
\hline \includegraphics[width=1.6cm]{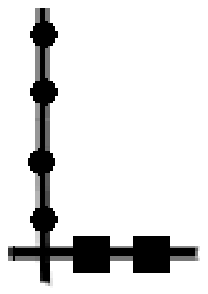} & $\Doii$ & $1$ \\
\hline \includegraphics[width=1.6cm]{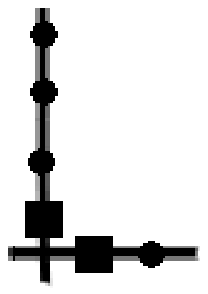} & $\Dor$ & $4$ \\
\hline
\end{tabular}
\end{center}

\begin{center}
\begin{tabular}{|m{4cm}|m{3cm}|m{3cm}|}
\hline Bound. Div. of $\bM_{0,[2,4]}$ & is mapped to& $M=N$ \\
\hline \hline \includegraphics[width=1.6cm]{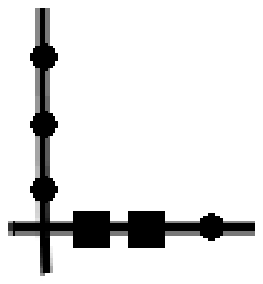} & $\Di$ & $6$ \\
\hline \includegraphics[width=1.6cm]{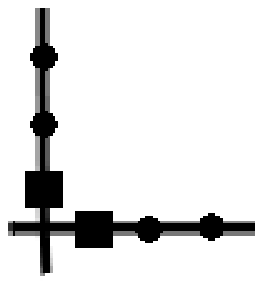} & $\Dii$ & $9$ \\
\hline
\end{tabular}
\end{center}

\begin{center}
\begin{tabular}{|m{4cm}|m{3cm}|m{3cm}|}
\hline Bound. Div. of $\bM_{0,[3,3]}$ & is mapped to& $M=N$ \\
\hline \hline \includegraphics[width=1.6cm]{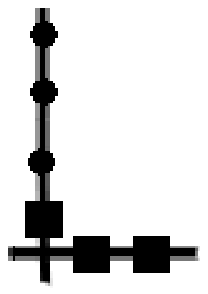} & $\Aop$ & $4$ \\
\hline \includegraphics[width=1.6cm]{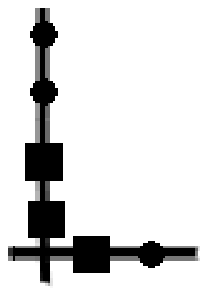} & $\Bop$ & $3$ \\
\hline
\end{tabular}
\end{center}

\begin{center}
\begin{tabular}{|m{4cm}|m{3cm}|m{3cm}|}
\hline Bound. Div. of $\bM_{0,[3,3]}$ & is mapped to& $M=N$ \\
\hline \hline\includegraphics[width=1.6cm]{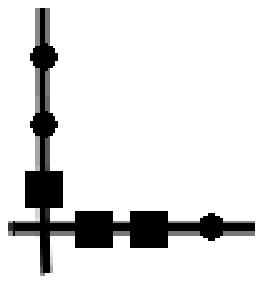} & $\Aip$ & $9$ \\
\hline \includegraphics[width=1.6cm]{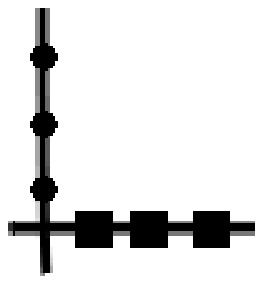} & $\Bip$ & $1$ \\
\hline
\end{tabular}
\end{center}

\begin{center}
\begin{tabular}{|m{4cm}|m{3cm}|m{3cm}|}
\hline Bound. Div. of $\bM_{0,[1,5]}$ & is mapped to& $M=N$ \\
\hline \hline \includegraphics[width=1.6cm]{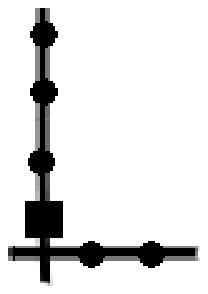} & $\Aom$ & $4$ \\
\hline \includegraphics[width=1.6cm]{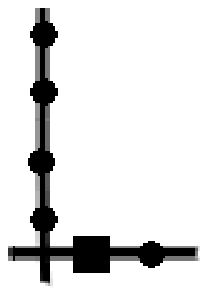} & $\Bom$ & $1$ \\
\hline \hline\includegraphics[width=1.6cm]{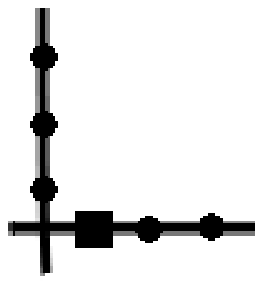} & $\Aip$ & $6$ \\
\hline
\end{tabular}
\end{center}

Now we can determine how $f_R$, $f_+$ and $f_-$ behave on the
boundary components of $\bM_{0,6}$. Using the notation introduced in
Definition \ref{short}, all these boundary components are of the
form $[i_1,i_2]$ or $[j_1,j_2,j_3]$ ($i_1,i_2,j_1,j_2,j_3 \in
\{1,2,3,4,5,6\}$). To which component a boundary component of
$\bM_{0,6}$ is mapped, can be seen using the tables above. The
degree of the map on a given boundary component one gets as in the
following example: The boundary component $[3,4]$ is mapped to
$\Doi$. A general point of $[3,4]$ is thus mapped by $f_R$ to a
point of $\Doi \subset \bR_2$ corresponding in $\bM_{0,[2,4]}$ to a
diagram of the form

\begin{center} \includegraphics[width=2cm]{doi.eps} \end{center}

One gets that the degree of $f_R$ on $[3,4]$ is $4$ by counting how
many nonisomorphic possibilities there are to assign indices
$1,...,6$ to the marked points of the diagram, such that the dots
get $3,4,5,6$, the squares get $1,2$ and such that $3$ and $4$ go to
the component with only two marked points. There are $8$
possibilities, but swapping $3$ and $4$ yields isomorphic objects.

Behavior of $f_R:\bM_{0,6} \overset{48:1}{\longrightarrow} \bR_2$. For arbitrary $b_1,b_2 \in \{3,4,5,6\}$ we have:
\begin{itemize}
  \item Boundary components of the form $[b_1,b_2]$ are mapped $4:1$ (each) onto $\Doi$.
  \item The boundary component $[1,2]$ is mapped $24:1$ onto $\Doii$.
  \item Boundary components of the form $[1,b_1]$ or $[2,b_1]$ are mapped $6:1$ (each) onto $\Dor$.
  \item Boundary components of the form $[1,2,b_1]$ are mapped $12:1$ (each) onto $\Di$.
  \item Boundary components of the form $[1,b_1,b_2]$ (or equivalently $[2,b_1,b_2]$) are mapped $8:1$ (each) onto $\Dii$.
\end{itemize}

Behavior of $f_+:\bM_{0,6} \overset{72:1}{\longrightarrow} \bS_2^+$. For arbitrary $a_1,a_2 \in \{1,2,3\}$ and $b_1,b_2 \in \{4,5,6 \}$ we have:
\begin{itemize}
  \item Boundary components of the form $[a_1,a_2]$ or $[b_1,b_2]$ are mapped $6:1$ (each) onto $A_0^+$.
  \item Boundary components of the form $[a_1,b_1]$ are mapped $8:1$ (each) onto $B_0^+$.
  \item Boundary components of the form $[a_1,a_2,b_1]$ (or equivalently $[a_1,b_1,b_2]$) are mapped $8:1$ (each) onto $A_1^+$.
  \item The boundary component $[1,2,3]$ is mapped $72:1$ onto $B_1^+$.
\end{itemize}

Behavior of $f_-:\bM_{0,6} \overset{120:1}{\longrightarrow} \bS_2^-$. For arbitrary $b_1,b_2 \in \{2,3,4,5,6\}$:
\begin{itemize}
  \item Boundary components of the form $[1,b_1]$ are mapped $24:1$ (each) onto $B_0^-$.
  \item Boundary components of the form $[b_1,b_2]$ are mapped $6:1$ (each) onto $A_0^-$.
  \item Boundary components of the form $[1, b_1,b_2]$ are mapped $12:1$ (each) onto $A_1^-$.
\end{itemize}

We now use this to compute:

\begin{Lemma} \label{a8}
There are the following relations between cycle classes on our
moduli spaces:

(i) In $A_{2,\mbb{Q}}(\bR_2)$: $\doi + 6 \doii -3 \dor +12 \di -8 \dii = 0$

(ii) In $A_{2,\mbb{Q}}(\bS_2^+)$: $3 \aop -4 \bop -8 \aip + 72 \bip =0$
\end{Lemma}

\textbf{Proof:} (i): Using equation (\ref{a9}) from Summary \ref{n2} with $i,j,k,l:=1,2,3,4$ we get
\[ [1,2]+[1,2,5]+[1,2,6]+[1,2,5,6] = [1,3]+[1,3,5]+[1,3,6]+[1,3,5,6]\]
which is the same as
\[ 0= [1,2]+[1,2,5]+[1,2,6]+[3,4]- [1,3]-[1,3,5]-[1,3,6]-[2,4]\]
Pushing this relation forward by $f_R$ we get:
\[ 0 = 24 [\Doii] + 12[\Di] + 12[\Di] + 4 [\Doi] - 6 [\Dor]- 8[\Dii] -8[\Dii]- 6 [\Dor]\]
\[ =4 [\Doi] + 24 [\Doii] - 12 [\Dor] + 24[\Di] -16[\Dii]\]
Using the automorphism numbers from the tables in the appendix, this
can be written as
\[0= 8 \doi + 48 \doii - 24 \dor + 96 \di -64 \dii\]
\[\Leftrightarrow  0= \doi + 6 \doii - 3 \dor + 12 \di -8 \dii \]

(ii): Using equation (\ref{a9}), this time with $i,j,k,l:=1,2,4,5$, we get
\[ [1,2]+[1,2,3]+[1,2,6]+[1,2,3,6] = [1,4]+[1,3,4]+[1,4,6]+[1,3,4,6]\]
Pushing this relation forward by $f_+$, and proceeding like in part (i) we get:
\[0= 24 \aop - 32 \bop -64 \aip +576 \bip\]
\[\Leftrightarrow  0= 3 \aop - 4 \bop -8 \aip +72 \bip \]
$\square$

\subsection{Morphisms to the boundary components of $\bR_2$ and $\bS_2$}

Now we come to several finite surjective morphisms from other moduli
spaces to different boundary components of $\bR_2$, $\bS_2^+$ and
$\bS_2^-$. Later they will be used to determine relations between
intersection products of boundary components via push-pull formula.

\subsubsection{Morphisms from $\bM_{0,5}$} \label{ss6}

Fist we define a Morphism $c: \bM_{0,5} \times \bM_{0,3} \rightarrow
[5,6] \subseteq \bM_{0,6}$. ($[5,6]$ is one of the boundary divisors
of $\bM_{0,6}$, c.f. Definition \ref{short}.) To the pair of
$[(C;(q_0,...,q_4))] \in \bM_{0,5}$ and $[(C';(q'_0,...,q'_2)] \in
\bM_{0,3}$ the morphism $c$ assigns $[D;(p_1,...,p_6)] \in [5,6]
\subset \bM_{0,6}$, where $D$ is the curve obtained from $C$ and
$C'$ by gluing the points $q_0$ and $q'_0$, and where $p_1,...,p_4$
are defined as the images of $q_1,...,q_4$ at $D$, and $p_5$ resp.
$p_6$ are defined as the images of $q'_1$ resp. $q'_2$. $\bM_{0,3}$
is just a point, so there is an isomorphism $i: \bM_{0,5}
\rightarrow \bM_{0,5} \times \bM_{0,3}$. The composed map $c \circ
i$ is a finite degree $1$ morphism onto $[5,6]$. We compose this
morphism with $f_R$ and get a finite Morphism:

\[ h_0': \bM_{0,5} \overset{4:1} \longrightarrow \Doi \]

$h_0'$ is $4:1$ because that is the degree of $f_R$ on $[5,6]$ (c.f.
section \ref{sec1}).

By composing $c \circ i$ with $f_-$ one gets a morphism

\[ h_0^{\alpha}: \bM_{0,5} \overset{6:1} \longrightarrow \Aom \]

Similar to what was done in section \ref{ss5} for $f_R$, $f_+$ and
$f_-$, one can determine the behavior of these two morphisms on the
boundary of $\bM_{0,5}$. We describe the images of the boundary
components in terms of the classes of closed strata of the
stratification by topological type of $\bR_2$ resp. $\bS_2^-$. These
strata are described in the appendix. The boundary divisors of
$\bM_{0,5}$ are (for our choice of the indices of the marked points)
all of the form $[i_1,i_2]$  ($i_1,i_2 \in \{0,1,2,3,4\}$).

Behavior of $h_0':\bM_{0,5} \overset{4:1}{\longrightarrow} \Doi \subset \bR_2$. For arbitrary $a \in \{1,2\}$ and  $b \in \{3,4 \}$:
\begin{itemize}
  \item The boundary component $[1,2]$ is mapped $2:1$ onto $E^{\prime , \prime \prime} = \Doi \cap \Doii$.
  \item Boundary components of the form $[a,b]$ are mapped $1:1$ (each) onto $E^{\prime , r} = \Doi \cap \Dor$.
  \item The boundary component $[3,4]$ is mapped $2:1$ onto $E^{\prime , \prime}$.
  \item Boundary components of the form $[0,a]$ are mapped $2:1$ (each) onto $F_{1:1}'= \Doi \cap \Dii$.
  \item Boundary components of the form $[0,b]$ are mapped $2:1$ (each) onto $F_1'= \Doi \cap \Di$.
\end{itemize}

Behavior of $h_0^{\alpha}:\bM_{0,5} \overset{6:1}{\longrightarrow} \Aom \subset \bS_2^-$. For arbitrary $b_1,b_2 \in \{2,3,4 \}$:
\begin{itemize}
  \item Boundary components of the form $[b_1,b_2]$ are mapped $2:1$ (each) onto $C^-$.
($2:1$ because two nonisomorphic diagrams of $\bM_{0,5}$ are assigned two different
but isomorphic diagrams of $\bM_{0,[1,5]} \cong \bS_2^-$.)
  \item Boundary components of the form $[1,b_1]$ are mapped $2:1$ (each) onto $D^- = \Aom \cap \Bom$.
  \item The boundary component $[0,1]$ is mapped $6:1$ onto $X^-$.
  \item Boundary components of the form $[0,b_1]$ are mapped $2:1$ (each) onto $Y^-$.
\end{itemize}

we use this to compute:

\begin{Lemma} \label{a10}
There are the following relations between classes in the Chow ring
of our moduli spaces:

(i) In $A_{1,\mbb{Q}}(\bR_2)$: $2 \doi \doii + 4 \doi \di - 4 \doi \dii - \doi \dor = 0$

(ii) In $A_{1,\mbb{Q}}(\bS_2^-)$: $16 [X^-]_Q + [C^-]_Q - 4 \aom
\aim - \aom \bom = 0$

(iii) In $A_{1,\mbb{Q}}(\bR_2)$: $[E^{\prime, r}]_Q = 2 [E^{\prime, \prime}]_Q +[E^{\prime, \prime \prime}]_Q$
\end{Lemma}

\textbf{Proof:} (i): Using equation \ref{a9} with $i,j,k,l:=0,1,2,3$ we get
\[ [0,1]+[2,3]= [0,3]+[1,2]\]
Pushing this relation forward by $h_0'$ we get:
\[ 0 = 2 [\Doi \cap \Di] + 2[\Doi \cap \Doii] - 2[\Doi \cap \Dii] - [\Doi \cap \Dor ]\]
Using the automorphism numbers from the appendix this can be written as
\[0=  8 \doi \di + 4 \doi \doii - 8 \doi \dii -2 \doi \dor\]
\[\Leftrightarrow  2 \doi \doii + 4 \doi \di - 4 \doi \dii - \doi \dor = 0 \]

(ii): We again use the equation
\[ 0= [0,3]+[1,2]-[0,1]-[2,3]\]
and now push it forward by $h_0^{\alpha}$. Then proceeding as above, we arive at  
\[ 0 = 12 [X^-]_Q + [C^-]_Q - 4[Y^-]_Q - \aom \bom\]
Now we use that $\Aom \cap \Aim = X^- \cup Y^-$ is a proper
intersection. We can treat all proper
intersections of $Q$-classes of closed strata of the stratifications
by topological type, as transversal intersections, as those closed
strata meet transversally on the deformation space (c.f. Summary
\ref{p2} (iv)). Thus $\aom \aim = [X^-]_Q + [Y^-]_Q$. Using this one
can rewrite the equation as
\[0 = 16 [X^-]_Q + [C^-]_Q - 4\aom \aim - \aom \bom\]

(iii) Using equation \ref{a9} with $i,j,k,l:=1,2,3,4$ we get
\[ [1,2]+[3,4]= [1,3]+[2,4]\]
Pushing this relation forward by $h_0'$ and using the automorphism numbers from the appendix we get:
\[ 4 [E^{\prime,\prime \prime}]_Q + 8[E^{\prime,\prime}]_Q = 2[E^{\prime,r} ]_Q\]
\[\Leftrightarrow  [E^{\prime,\prime \prime}]_Q + 2[E^{\prime,\prime}]_Q = [E^{\prime,r} ]_Q \]
$\square$

\subsubsection{Other morphisms to the boundary components} \label{ss7}

For $\bR_2$ we will use the following morphisms. We describe how they behave on general points.

\[ \ti: \bM_{1,1} \times \bR_{1,1} \overset{1:1}{\longrightarrow} \Di \]
For $[(X;p)] \in \bM_{1,1}$ and $[(Y;\mc{L};b;q)] \in \bR_{1,1}$ the
image of the pair is the point in $\Di$ parametrizing the following
Prym curve $(X'; \mc{L}')$: The quasistable curve $X'$ is generated
by gluing the points $p$ and $q$ on the curves $X$ and $Y$. The Prym
sheaf $\mc{L}'$ is obtained from the trivial sheaf on $X$ and the
Prym sheaf $\mc{L}$ on $Y$, by identifying the fibers over $p$ resp.
$q$. All possible choices of identification yield isomorphic Prym
sheaves.

\[\tii: \bR_{1,1} \times \bR_{1,1} \overset{2:1}{\longrightarrow} \Dii \]
This morphism is defined analogously to $\ti$ . It is of degree $2$
since a pair $([(X;\mc{L};b;p)], [(X';\mc{L}';b';p')]) \in \bR_{1,1}
\times \bR_{1,1}$
and the transposed pair are mapped to the same point in $\Dii$. \\

\[\toii: \bM_{1,2} \overset{1:1}{\longrightarrow} \Doii \]
A point $[(X; p,q)] \in \bM_{1,2}$ is mapped to the point
parameterizing the following Prym curve $(X'; \mc{L})$: The
underlying quasistable curve $X'$ is obtained by gluing the points
$p$ and $q$. There are two ways to glue the fibers of the trivial
bundle of $X$ over the points $p$ and $q$ such that a Prym bundle on
$X'$ is obtained. One way yields the trivial bundle on $X'$, the
other one yields the nontrivial Prym bundle $\mc{L}$.

\[\tor: \bR_{1,2}^{(-1,-1)} \overset{1:1}{\longrightarrow} \Dor \]
A point $[(X; \mc{L}; p,q)] \in \bR_{1,2}^{(-1,-1)}$ is mapped to
the point parametrizing the following Prym curve $(X'; \mc{L}')$:
The underlying quasistable curve $X'$ is obtained by gluing the
points $p$ and $q$, and then blowing up the node. $\mc{L}'$ is the
Prym bundle on $X$, such that if $\tilde{X}$ is stable subcurve of
$X$ and $E$ the exceptional component, $\mc{L}'_{|\tilde{X}} \cong
\mc{L}$ and $\mc{L}'_{|E} \cong \GO_E(1)$.

\[\toi: \bM_{0,(2,2,1)} \overset{1:1}{\longrightarrow} \Doi \]
The morphism $h_0'$ factors through $\bM_{0,(2,2,1)}$, and we use this to define $\toi$.

For $\bS_2^+$ we will use the following morphisms.

\[ \roa: \bS_{1,2}^{(1,1)} \overset{1:1}{\longrightarrow} \Aop \]
A point $[(X; \mc{L}; b; p,q)] \in \bS_{1,2}^{1,1}$ is mapped to the
point parametrizing  the following spin curve $(X'; \mc{L}')$: The
underlying quasistable curve $X'$ is obtained by gluing the points
$p$ and $q$. There are two ways to glue the fibers of the the bundle
$\mc{L}$ of $X$ over the points $p$ and $q$ such that a spin bundle
on $X'$ is obtained. One way yields an odd bundle, the other one the
even bundle $\mc{L}'$. (This is implicit in \cite{MR1082361},
Example 3.2)

\[ \rob: \bS_{1,2}^+ \overset{1:1}{\longrightarrow} \Bop \]
Defined analogously to $\tor$.

\[ \ria: \bS_{1,1}^+ \times \bS_{1,1}^+ \overset{2:1}{\longrightarrow} \Aip \]
Defined analogously to $\ti$, but the node is blown up.

\[ \rib: \bS_{1,1}^- \times \bS_{1,1}^- \overset{2:1}{\longrightarrow} \Bip \]
Defined analogously to $\ria$

For $\bS_2^-$ there are the following morphisms.

\[ \hoa: \bS_{1,2}^{(1,1)} \overset{1:1}{\longrightarrow} \Aom \]
Defined analogously to $\roa$.

\[ \hob: \bS_{1,2}^- \overset{1:1}{\longrightarrow} \Bom \]
Defined analogously to $\rob$.

\[ \hia: \bS_{1,1}^+ \times \bS_{1,1}^-  \overset{1:1}{\longrightarrow} \Aim \]
Defined analogously to $\ria$.

Now we gather facts about some of the moduli spaces of pointed
curves that the domains of the morphisms just defined consist of.
Especially this will be facts about the rational Chow groups of
these spaces.

\begin{enumerate}
\item $\bM_{1,1}$ has only one boundary divisor: $\tilde{\Delta}_0$. It parametrizes curves with one node.
The corresponding $Q$-class we call $\tilde{\delta_0}:= [\tilde{\Delta}_0]_Q$.

\item $\bR_{1,1}$ has boundary divisors $\tilde{\Doii}$ and $\tilde{\Dor}$, defined analogously to $\Doii$ and $\Dor$.
The corresponding $Q$-classes we call $\tilde{\doii}$ and $\tilde{\dor}$. $\bR_{1,1}$ is isomorphic to $\mbb{P}^1$, thus $\tilde{\doii} = \tilde{\dor}$ in the Chow group.

\item $\bM_{1,2}$ has boundary divisors $\hat{\Delta}_0$ and $\hat{\Delta}_1$. A curve parametrized by a general point of $\hat{\Delta}_0$ is irreducible with one node.
 A general curve parametrized by $\hat{\Delta}_1$ consists of two irreducible components, one smooth elliptic curve and one smooth rational curve with two marked points.
The corresponding $Q$-classes we call $\hat{\delta_0}$ and $\hat{\delta_1}$.

\item $\bR_{1,2}$ has  boundary divisors $\hat{\Doii}$, $\hat{\Dor}$ and $\hat{\Di}$.
Where $\hat{\Doii}$ and $\hat{\Dor}$ are defined analogously to
$\Doii$ and $\Dor$. For a Prym curve $(X;\mc{L}; b;p,q)$
parametrized by a general point of $\hat{\Di}$, $X$ consists of two
irreducible components, one smooth elliptic curve and one smooth
rational curve with two marked points. The Prym sheaf $\mc{L}$ is
nontrivial restricted to the elliptic curve and (necessarily)
trivial restricted to the rational curve. The $Q$-classes
$\hat{\doii}$ and $\hat{\dor}$ are equivalent in the Chow group,
because they are the pullbacks of the corresponding classes on
$\bR_{1,1}$.

\item $\bS_{1,1}^-$ and $\bS_{1,2}^-$ are just $\bM_{1,1}$ respectively $\bM_{1,2}$ because a odd Prym sheaf on a genus $1$ curve is trivial.
\textbf{In later computations, we will usually replace $\bS_{1,1}^-$
and $\bS_{1,2}^-$ by $\bM_{1,1}$ respectively $\bM_{1,2}$ without
further mentioning}.

\item $\bS_{1,1}^+$: The boundary divisors are $\tilde{A}_0^+$ and $\tilde{B}_0^+$. Defined analogously to $\Aop$ and $\Bop$.
The corresponding $Q$-classes $\tilde{\alpha}_0^+$ and $\tilde{\beta}_0^+$ are equivalent in the Chow group, since
$\bS_{1,1}^+ \cong \mbb{P}^1$.

\item $\bS_{1,2}^+$: The boundary divisors are $\hat{A}_0^+$, $\hat{B}_0^+$ and $\hat{A}_1^+$. The $Q$-classes $\hat{\alpha}_0^+$  and $\hat{\beta}_0^+$ are equivalent in the Chow group, since they are the pullbacks of the corresponding classes on $\bS_{1,1}$.

\item $\bS_{1,2}^{(1,1)}$: There are, among others, the boundary divisors $\check{A}_0$ and $\check{B}_0$ whose general points parametrize irreducible curves with one node that is blown up in the case of $\check{B}_0$. The $Q$ classes $\check{\alpha}_0$ and $\check{\beta}_0$ are not equivalent.
\end{enumerate}

The facts just listed are probably all known (for some of them c.f.
\cite{MR2526309}, Page 8, and \cite{MR2560968}). One way of proving
them, is to use that the moduli spaces of curves with one marked
points appearing, are all isomorphic to certain quotients of
$\bM_{0,4}$. The moduli spaces of curves with two marked points
appearing, are, after forgetting the order of the two marked points,
isomorphic to certain quotients of $\bM_{0,5}$. For an example look at
Part (ii) of the following Lemma. Forgetting the order of the two
marked points on the genus $1$ curves does not change the coarse
moduli spaces.

\begin{Lemma} \label{x1}
 (i) Let
 \[\pi_{(2,2,1)}: \bM_{0,5} \rightarrow \bM_{0,(2,2,1)}, \quad [(X;(p_1,....,p_4,p_0))] \mapsto [(X;(\{p_1,p_2\},\{p_3,p_4\},\{p_0\}))] \]
 be the quotient morphisms of Definition \ref{rd1}, and let $a \in \{1,2\}$ and $b \in \{3,4\}$ be arbitrary. We define
 \[ C'':= \pi_{(2,2,1)}([1,2]), \quad C':= \pi_{(2,2,1)}([3,4]), \quad C^r:=
\pi_{(2,2,1)}([a,b]),\]
\[ C_{1:1}:= \pi_{(2,2,1)}([a,0]), \quad C_1:=
\pi_{(2,2,1)}([b,0]) \]
 These images are independent of the choice of $a$ and $b$, which implies that the moduli space $\bM_{0,(2,2,1)}$
 has exactly the five
 boundary components $C'$, $C''$, $C^r$, $C_1$  and $C_{1:1}$.

(ii) There is an isomorphism $\bM_{0,[4,1]} \rightarrow \bM_{1,[2]}
\cong \bM_{1,2}$. By combining this with the forgetful morphism
$\bM_{0,(2,2,1)} \rightarrow \bM_{0,[4,1]}$ we define a finite
surjective morphism $\theta: \bM_{0,(2,2,1)} \rightarrow \bM_{1,2}$
\end{Lemma}

\textbf{Proof:} (i): Easy to check. (For the notation used, c.f.
Definition \ref{short}.)

(ii): To a point $[(D; A,p)] \in \bM_{0,[4,1]}$, let $f:Y
\rightarrow D$ be the admissible $2:1$ cover of $(D;A)$, and let $Q$
be the set $f^{-1}(p)$. Then $[(D; A,p)] \mapsto [(Y; Q)]$ defines a
morphism $\theta': \bM_{0,[4,1]} \rightarrow \bM_{1,[2]} \cong
\bM_{1,2} $. It is easy to check that it is 1:1 on the locus of
smooth curves. Since both moduli spaces are normal projective
varieties this suffices to prove that $\theta'$ is an isomorphism.
$\square$

\begin{Lemma}\label{a11}
The following table shows the pushforwards of several classes by the morphisms defined in this section.

\begin{center}
\begin{tabular}{|c|c|c|}
\hline \textbf{Morphism} & \textbf{class} & \textbf{Pushforward}\\
\hline
\hline $\toi$ & $1$ & $2 \doi$\\
\hline $\toi$ & $c'$ & $2 [E^{\prime,\prime}]_Q$\\
\hline $\toi$ & $c''$ & $ \doi \doii$\\
\hline $\toi$ & $c^r$ & $2 \doi \dor$\\
\hline $\toi$ & $c_1$ & $4 \doi \di$\\
\hline $\toi$ & $c_{1:1}$ & $4 \doi \dii$\\

\hline $\toii$ & $1$ & $2 \doii$\\
\hline $\toii$ & $\hat \delta_0  $ & $2 \doi \doii$\\
\hline $\toii$ & $\hat \delta_1  $ & $2 \doi \di$\\

\hline $\ti$ & $\tilde{d}_0'' \otimes 1$ & $\doii \di$\\
\hline $\ti$ & $\tilde{d}_0'' \otimes 1$ & $\doii \di$\\
\hline $\ti$ & $\tilde{d}_0^r \otimes 1$ & $\dor \di$\\
\hline $\ti$ & $1 \otimes \tilde{\delta}_0$ & $\doi \di$\\

\hline $\tii$ & $\tilde{d}_0'' \otimes 1$ & $\doi \dii$\\
\hline $\tii$ & $1 \otimes \tilde{d}_0''$ & $\doi \dii$\\
\hline $\tii$ & $\tilde{d}_0^r \otimes 1$ & $\dor \dii$\\
\hline $\tii$ & $1 \otimes \tilde{d}_0^r$ & $\dor \dii$\\

\hline
\end{tabular}

\begin{tabular}{|c|c|c|}
\hline \textbf{Morphism} & \textbf{class} & \textbf{Pushforward}\\
\hline

\hline $\roa$ & $1$ & $2 \aop$\\
\hline $\roa$ & $\check{\alpha}_0$ & $4 [C^+]_Q$\\
\hline $\roa$ & $\check{\beta}_0$ & $2 \aop \bop$\\

\hline $\rob$ & $1$ & $2 \bop$\\
\hline $\rob$ & $\hat{\alpha}_0^+$ & $2 \aop \bop$\\
\hline $\rob$ & $\hat{\beta}_0^+$ & $4 [E]_Q$\\

\hline $\ria$ & $\tilde{\alpha}_0^+ \otimes 1$ & $2 \aop \aip$\\
\hline $\ria$ & $1 \otimes \tilde{\alpha}_0^+$ & $2 \aop \aip$\\
\hline $\ria$ & $\tilde{\beta}_0^+ \otimes 1$ & $2 \bop \aip$\\
\hline $\ria$ & $1 \otimes \tilde{\beta}_0^+$ & $2 \bop \aip$\\

\hline $\rib$ & $\tilde{\delta}_0 \otimes 1$ & $2 \aop \bip$\\
\hline $\rib$ & $1 \otimes \tilde{\delta}_0$ & $2 \aop \bip$\\

\hline
\end{tabular}

\begin{tabular}{|c|c|c|}
\hline \textbf{Morphism} & \textbf{class} & \textbf{Pushforward}\\
\hline
\hline $\hoa$ & $1$ & $2 \aom$\\
\hline $\hoa$ & $\check{\alpha}_0$ & $4 [C^-]_Q$\\
\hline $\hoa$ & $\check{\beta}_0$ & $2 \aom \bom$\\

\hline $\hob$ & $1$ & $2 \bom$\\
\hline $\hob$ & $\hat{\delta}_0$ & $2 \aop \bop$\\

\hline $\hia$ & $\tilde{\alpha}_0^+ \otimes 1$ & $2 [X^-]_Q$\\
\hline $\hia$ & $\tilde{\beta}_0^+ \otimes 1$ & $2 \bom \aim$\\
\hline $\hia$ & $1 \otimes \tilde{\delta}_0$ & $2 [Y^-]_Q$\\

\hline
\end{tabular}

\end{center}
\end{Lemma}

\textbf{Proof:} By counting the degree of the given morphism on the
given cycle, and comparing the automorphism number of an object
parametrized by a general point of the cycle, with the automorphism
number of the object parametrized by the image of such a point,
under the given morphism. $\square$

\subsection{Hodge classes}

Another type of cycle classes used in our computation, beside
classes of closed strata of the stratifications according to
topological type, are first Chern classes of the Hodge bundles on
moduli spaces, and their pullbacks.

\begin{Def} \label{d}
 Let $\tilde{\pi}_R: \bR_{1,1} \longrightarrow  \bM_{1,1}$, $\tilde{\pi}^+: \bS_{1,1}^+ \longrightarrow  \bM_{1,1}$, $\hat{\pi}^+: \bS_{1,2}^+ \longrightarrow  \bM_{1,2}$,
and $\check{\pi}: \bS_{1,2}^{(1,1)} \longrightarrow  \bM_{1,2}$  be
the usual forgetful morphisms, and let $\theta: \bM_{0,(2,2,1)} \rightarrow \bM_{1,2}$ be the morphism of
Lemma \ref{x1} (ii). Let $\lambda$, $\tilde{\lambda}$ resp.
$\hat{\lambda}$  be the first Chern class of
the Hodge bundle on $\bM_2$, $\bM_{1,1}$ resp. $\bM_{1,2}$. 

We define classes :
\[ l:= (\pi_R)^{\ast} \lambda, \quad l^+:= (\pi_+)^{\ast} \lambda, \quad l^-:= (\pi_-)^{\ast} \lambda, \quad \tilde{l}:= (\tilde{\pi}_R)^{\ast} \tilde{\lambda}, \]
\[ \tilde{l}^+:= (\tilde{\pi}^+)^{\ast} \tilde{\lambda}, \quad \hat{l}^+:= (\hat{\pi}^+)^{\ast} \hat{\lambda}, \quad \check{l}:= (\check{\pi})^{\ast} \hat{\lambda}, \quad  \bar{l}:= \theta^* \hat{\lambda}\]
\end{Def}

\begin{Lemma} \label{7}
 We can describe the pullbacks of $l$, $l^+$ and $l^-$ by the boundary morphisms in the following way

 (i) $ (\ti)^{\ast} l = \tilde{\lambda} \otimes 1 + 1 \otimes \tilde{l}$

 (ii) $ (\tii)^{\ast} l = \tilde{l} \otimes 1 + 1 \otimes \tilde{l}$

 (iii) $ (\toi)^* l = \bar{l}$

 (iv) $ (\toii)^{\ast} l = \hat{\lambda} $

 (v) $ (\roa)^{\ast} l^+ = \check{l}$

 (vi) $ (\rob)^{\ast} l^+ = \hat{l}^+$

 (vii) $ (\ria)^{\ast} l^+ = \tilde{l}^+ \otimes 1 + 1 \otimes \tilde{l}^+$

 (viii) $ (\rib)^{\ast} l^+ = \tilde{\lambda} \otimes 1 + 1 \otimes \tilde{\lambda}$

 (ix) $ (\hoa)^{\ast} l^- = \check{l}$

 (x) $ (\hob)^{\ast} l^- = \hat{\lambda}$

 (xi) $ (\hia)^{\ast} l^- = \tilde{\lambda} \otimes 1 + 1 \otimes \tilde{l}^+$

\end{Lemma}
{}
\textbf{Proof:} First consider the commutative diagram
\[ \begin{xy}
     \xymatrix{ \bS_{1,2}^{(1,1)} \ar[d]_{\check{\pi}}  \ar[r]^{\roa} & \bS_2^+ \ar[d]^{\pi_+}
            \\ \bM_{1,2} \ar[r]^f & \bM_2}
   \end{xy} \]
where $f$ is the morphism corresponding to gluing the two marked
points on a curve. Because of the way $l^+$ and $\check{l}$ are
defined, it suffices to show $\hat{\lambda} = f^* \lambda$ in order
to prove (v). That this equation indeed is true, is shown in
\cite{MR717614},  § 10. The assertions (iii), (iv), (vi), (ix) and (x) can
be proved in the same way.

Now we consider the commutative diagram:

\[ \begin{xy}
     \xymatrix{ \bR_{1,1} \times \bR_{1,1} \ar[d]_{\tilde{\pi}_R \times \tilde{\pi}_R } \ar[r]^{\tii} & \bR_2 \ar[d]^{\pi_R}
            \\ \bM_{1,1} \times \bM_{1,1} \ar[r]^g & \bM_2}
   \end{xy} \]
Where $g$ is the morphism corresponding to gluing two genus 1
curves, each with one marked point, together at those marked points.
In [\cite{MR717614},  § 10. $g^* \lambda = \tilde{\lambda} \otimes 1
+ 1 \otimes \tilde{\lambda}$ is proven (there the notation is
slightly different). From this (i) follows. (ii), (vii), (viii) and
(xi) can be proved analogously.

$\square$

If $\lambda$ is the fist Chern class of the Hodge bundle on a
$\bM_{1,n}$, $n \ge 1$ arbitrary, then for $\delta_0$ the $Q$ class
of the divisor of $\bM_{1,n}$ parametrizing irreducible curves with
one node, $\lambda = \frac{1}{12} \delta_0$ (c.f. \cite{MR2526309}
Page 8). By pulling these relations back one obtains the following
equation:

\begin{Lemma} \label{8}
 $\quad$

 (i) $\tilde{\lambda} = \frac{1}{12} \tilde{\delta_0}$

 (ii) $\hat{\lambda} = \frac{1}{12} \hat{\delta_0}$

 (iii) $\bar{l} = \frac{1}{12} (2 c'+ 2c''+ 2c^r)$

 (v) $\tilde{l} = \frac{1}{12} \tilde{\theta}^{\ast} \tilde{\delta_0}= \frac{1}{12} (\tilde{\doii} + 2 \tilde{\dor}) = \frac{1}{4}  \tilde{\dor}$

 (vi) $\check{l} = \frac{1}{12} (\check{\theta})^{\ast} \hat{\delta_0} = \frac{1}{12}( \check{\alpha}_0 + 2 \check{\beta}_0)$

 (v) $\tilde{l}^+ = \frac{1}{12} \tilde{\theta}^{\ast} \tilde{\delta_0}= \frac{1}{12}( \tilde{\alpha}_0^+ + 2 \tilde{\alpha}_0^+)
 = \frac{1}{4} \tilde{\alpha}_0^+$

 (v) $\hat{l}^+ = \frac{1}{12}( \hat{\alpha}_0^+ + 2 \hat{\alpha}_0^+)
 = \frac{1}{4} \hat{\alpha}_0^+$

\end{Lemma}

\begin{Lemma} \label{a12}
All the following products are equal to $0$ in the rational Chow
rings they are contained in.

\[l^2 \doi, \quad l^2 \doii, \quad l^2 \dor, \quad (l^+)^2 \aop, \quad (l^+)^2 \bop, \quad (l^-)^2 \aom, \quad (l^-)^2 \bom \]
\end{Lemma}

\textbf{Proof:} Take for example $(l^+)^2 \aop$. Using the boundary
morphism $\roa: \bS_{1,2}^{(1,1)} \overset{1:1}{\longrightarrow}
\Aop$ and the fact that $\aop= \frac{1}{2} (\roa)_*(1)$ we can write
$(l^+)^2 \aop$ by the projection formula as $\frac{1}{2}
(\roa)_*(\roa)^* (l^+)^2$. According to Lemma \ref{7} $(\roa)^*
(l^+)=\check{l}$, thus $(l^+)^2 \aop = \frac{1}{2} (\roa)_*
(\check{l})^2$. By definition $\check{l}=(\check{\pi})^{\ast}
\hat{\lambda}$. But $\hat{\lambda}$ is, as shown in \cite{MR717614}
§ 10., equal to the pullback of $\tilde{\lambda}$ from $\bM_{1,1}$
to $\bM_{1,2}$. $\bM_{1,1}$ is one dimensional, thus
$(\tilde{\lambda})^2=0$. This implies $(\check{l})^2=0$, which
pushed forward by $\roa$ yields $(l^+)^2 \aop=0$. That the other
products listed in the Lemma are equal to $0$ can be proved
analogously. $\square$

\section{Computation of the rational cohomology} \label{s3}

\subsection{The rational Picard group} \label{ss8}

\begin{Lemma} \label{a4}

The Chow groups $A_{2,\mbb{Q}}(\bR_2)$, $A_{2,\mbb{Q}}(\bS_2^+)$ and
$A_{2,\mbb{Q}}(\bS_2^-)$, are isomorphic to the rational Picard
groups $Pic_{\mbb{Q}}(\bR_2)$, $Pic_{\mbb{Q}}(\bS_2^+)$ respectively
$Pic_{\mbb{Q}}(\bS_2^-)$, an they are  generated by the boundary
divisors of the moduli spaces. Furthermore the linear relations of
Lemma \ref{a8} are the only ones. Thus:

(i) $A_{2,\mbb{Q}}(\bR_2)= (\doi \mbb{Q} \oplus \doii \mbb{Q} \oplus \dor \mbb{Q}\oplus \di \mbb{Q} \oplus \dii \mbb{Q})/(\doi + 6 \doii -3 \dor +12 \di -8 \dii) \mbb{Q}$

(ii) $A_{2,\mbb{Q}}(\bS_2^+)= (\aop \mbb{Q} \oplus \bop \mbb{Q} \oplus \aip \mbb{Q} \oplus \bip \mbb{Q})/(3 \aop -4 \bop -8 \aip + 72 \bip) \mbb{Q}$

(iii) $A_{2,\mbb{Q}}(\bS_2^-)= \aom \mbb{Q} \oplus \bom \mbb{Q} \oplus \aim \mbb{Q}$
\end{Lemma}
{}
\textbf{Proof:} That the second Chow groups are generated by boundary divisors and are isomorphic to the rational Picard groups
is a special case of Corollary \ref{n3} (iv) resp. (iii).

It remains to show that there are no linear relations between the boundary classes other than those of lemma \ref{a8}.

To do this we compute the intersection numbers of all boundary
classes with the $Q$-classes of all the $2$-dimensional closed
strata of the stratifications of our moduli spaces according to
topological type. These are the cycles lying above the
cycles $\Delta_{00}$ and $\Delta_{01}$ of $\bM_2$ with respect to
the forgetful morphisms. They are described in the appendix.
For a codimension $1$ cycle $d$ and a codimension $2$ cycle $e$ we
take the intersection number to be the number $n$ such that
$de=n[x]$ where $x$ is a general point of the moduli space. Note
that in the definition we use the class $[x]$, not $[x]_Q$, to be
consistent with \cite{MR717614}. For $\bR_2$ we get the intersection
numbers:

\[\begin{tabular} {|c|c||c|c|c|c|c|}
 \hline
 Underlying stratum of $\bM_2$ & stratum class & $\doi$ & $\doii$ & $\dor$ & $\di$ & $\dii$ \\
 \hline \hline  $\Delta_{00}$ & $[E^{\prime,\prime}]_{Q}$ & $-\frac{1}{2}$ & $\frac{1}{4}$ & $0$ & $0 $ & $\frac{1}{8}$\\
 \hline  $\Delta_{00}$ & $[E^{\prime,\prime \prime}]_{Q}$ & $0$ & $-\frac{1}{2}$ & $0$ & $\frac{1}{4}$ & $0$\\
 \hline  $\Delta_{00}$ & $[E^{\prime,r}]_{Q}$ & $-1$ & $0$ & $0$ & $\frac{1}{4}$ & $\frac{1}{4}$\\
 \hline  $\Delta_{00}$ & $[E^{r,r}]_{Q}$ & $\frac{1}{4}$ & $0$ & $-\frac{1}{4}$ & $0$ & $\frac{1}{8}$\\

 \hline \hline  $\Delta_{01}$ & $[F_1']_{Q}$ & $0$ & $\frac{1}{4}$ & $\frac{1}{4}$ & $-\frac{3}{48}$ & $0$\\
 \hline  $\Delta_{01}$ & $[F_1'']_{Q}$ & $\frac{1}{4}$ & $0$ & $0$ & $-\frac{1}{48}$ & $0$\\
 \hline  $\Delta_{01}$ & $[F_1^r]_{Q}$ & $\frac{1}{4}$ & $0$ & $0$ & $-\frac{1}{48}$ & $0$\\
 \hline  $\Delta_{01}$ & $[F_{1:1}']_{Q}$ & $\frac{1}{4}$ & $0$ & $\frac{1}{4}$ & $0$ & $-\frac{3}{48}$\\
 \hline  $\Delta_{01}$ & $[F_{1:1}^r]_{Q}$ & $\frac{1}{4}$ & $0$ & $\frac{1}{4}$ & $0$ & $-\frac{3}{48}$\\
 \hline
\end{tabular}\]

If we have a linear relation $\alpha_1 \doi + \alpha_2 \doii +
\alpha_3 \dor + \alpha_4 \di +\alpha_5 \dii = 0$ between the boundary
components, the vector $\alpha = (\alpha_1,...,\alpha_5)$ has to lie
in the kernel of the $9 \times 5$ matrix formed by the intersection
numbers in the table above. One can check, that this matrix has rank
$4$ and thus has $1$-dimensional kernel, and that the relation $\doi
+ 6 \doii -3 \dor +12 \di -8 \dii$ indeed lies in its kernel.

For $\bS^+$ the intersection numbers are:

\[\begin{tabular} {|c|c||c|c|c|c|}
 \hline
 Underlying stratum of $\bM_2$ & stratum class & $\aop$ & $\bop$ & $\aip$ & $\bip$ \\
 \hline \hline  $\Delta_{00}$ & $[C^+]_{Q}$ & $-1$ & $\frac{1}{4}$ & $\frac{1}{16}$ & $\frac{1}{16}$\\
 \hline  $\Delta_{00}$ & $[D^+]_{Q}$ & $0$ & $-\frac{1}{4}$ & $\frac{1}{8}$ & $0$\\
 \hline  $\Delta_{00}$ & $[E]_{Q}$ & $0$ & $-\frac{1}{8}$ & $\frac{1}{16}$ & $0$\\

 \hline \hline  $\Delta_{01}$ & $[X^+]_{Q}$ & $\frac{1}{8}$ & $\frac{1}{8}$ & $-\frac{3}{192}$ & $0$\\
 \hline  $\Delta_{01}$ & $[Y^+]_{Q}$ & $\frac{1}{8}$ & $0$ & $0$ & $-\frac{1}{192}$\\
 \hline  $\Delta_{01}$ & $[Z^+]_{Q}$ & $\frac{1}{8}$ & $\frac{1}{8}$ & $-\frac{3}{192}$ & $0$\\
 \hline
\end{tabular}\]

One can check that the $6 \times 4$ matrix formed by the intersection numbers, has rank $3$, and that
$3 \aop -4 \bop -8 \aip + 72 \bip$ lies inside the kernel.

For $\bS^-$ the intersection numbers are:

\[\begin{tabular} {|c|c||c|c|c|}
 \hline
 Underlying stratum of $\bM_2$ & stratum class & $\aom$ & $\bom$ & $\aim$ \\
 \hline \hline  $\Delta_{00}$ & $[C^-]_{Q}$ & $-1$ & $\frac{1}{4}$ & $\frac{1}{8}$ \\
 \hline  $\Delta_{00}$ & $[D^-]_{Q}$ & $0$ & $-\frac{1}{4}$ & $\frac{1}{8}$ \\

 \hline \hline  $\Delta_{01}$ & $[X^-]_{Q}$ & $\frac{1}{8}$ & $0$ & $-\frac{1}{192}$\\
 \hline  $\Delta_{01}$ & $[Y^-]_{Q}$ & $\frac{1}{8}$ & $\frac{1}{8}$ & $-\frac{3}{192}$ \\
 \hline  $\Delta_{01}$ & $[Z^-]_{Q}$ & $\frac{1}{8}$ & $0$ & $-\frac{1}{192}$\\
 \hline
\end{tabular}\]

The $5 \times 3$ matrix formed by the intersection numbers has rank $3$.

As examples we will compute some intersection numbers from the
tables above. The other numbers can be computed analogously. From
\cite{MR717614}, Theorem $10.1$, we know that $\delta_0
[\Delta_{00}]_Q =- \frac{1}{4} p$, $\delta_1  [\Delta_{00}]_Q =
\frac{1}{8} p$, $\delta_1 [\Delta_{01}]_Q  = -\frac{1}{48} p $ and
$\delta_0 [\Delta_{01}]_Q  = \frac{1}{4} p$, where $p$ is the class
$[y]$ of a general point of $\bM_2$.

For $\overline{X} \in \{ \bR_2, \bS_2^+, \bS_2^- \}$ let $S$ be one
of the codimension $2$ cycles on $\overline{X}$ listed in the tables
above. If $\pi: \overline{X} \rightarrow \bM_2$ is the forgetful
morphism, then $\pi_*  S = m D$ for some $m \in \mbb{Q}$, and for
$D$ the $Q$-class of the reduced image of $S$ under $\pi$, thus $D =
[\Delta_{00}]_Q$ or $D= [\Delta_{01}]_Q$. The number $m$ is listed
for all cycles $S$ in the appendix. Thus one can compute the
intersection number $n$ of $S$ with the pullback of $\delta_i$
($i=0,1$) by using the forgetful map $\pi$ and the projection
formula:
\[ \pi^* \delta_i S  = n [x] \quad \Leftrightarrow \quad \delta_i \pi_* S = n [y] = n p\]
\[ \Leftrightarrow \quad  m D \delta_i = n p\]
Where $D \delta_i$ is one of the four known intersections on $\bM_2$ mentioned above.

For the example $E^{\prime,\prime}$ we have $(\pi_R)_* [E^{\prime,\prime}]_Q =  [\Delta_{00}]_Q$, thus
$\pi^* \delta _0  [E^{\prime,\prime}]_Q = - \frac{1}{4} [x]$ and $\pi^* \delta _1  [E^{\prime,\prime}]_Q = \frac{1}{8} [x]$.

We also have $\Dor \cap E^{\prime,\prime} = \Di \cap
E^{\prime,\prime} = \emptyset$ (as one can show using the
description of these strata in the appendix), so the corresponding
intersection numbers are $0$.
 Using $(\pi_R)^* \delta_0 = \doi + \doii+ 2 \dor$ and
$(\pi_R)^* \delta_1 = \di + \dii$, we get $\di [E^{\prime,\prime}]_Q = \frac{1}{8} [x]$ and
\begin{equation} \label{p}
(\doi + \doii)  [E^{\prime,\prime}]_Q = - \frac{1}{4} [x]
\end{equation}

The intersection $\Doii \cap E^{\prime,\prime} = G'$ (use
description in the appendix) is proper, so by Summary \ref{p2} (iv)
we can treat the intersection as transversal and we get $\doii
[E^{\prime,\prime}]_Q = [G']_Q$. $G'$ consist of one point, and the
corresponding Prym curve has $4$ automorphisms (c.f. appendix), thus
$\doii [E^{\prime,\prime}]_Q = \frac{1}{4} [x]$. By plugging this
into equation (\ref{p}) we obtain the last intersection number $\doi
[E^{\prime,\prime}]_Q = - \frac{1}{2} [x]$.

All rows in the above tables can be computed in this way, except for
the ones containing the intersection numbers of $E^{\prime, \prime
\prime}$, $E^{\prime, r}$ and $D^-$. In computing the first two one has to use
additionally the relation $[E^{\prime, r}]_Q = 2 [E^{\prime,
\prime}]_Q +[E^{\prime, \prime \prime}]_Q$. For the intersections with $[D^-]_Q$ one uses the relation $12[X^-]_Q+[C^-]_Q-4[Y^-]_Q=[D^-]_Q$. Both relations are proven
in Lemma \ref{a10}. $\square$

\textbf{Remark:} In \cite{MR2526309}, Page 5-6, it is claimed
that the boundary components of $S_2^+$ (and $S_2^-$) are
independent, which results in wrong Betti (and Hodge) numbers computed
for $S_2^+$. It is claimed that Cornalba's proof of independence of the boundary
classes for genus $g \ge 3$ in \cite{MR1082361}, can also be applied
to $g=2$. Cornalba's proof works similar to the proof of the lemma
above by computing intersections of the boundary classes with
various test curves. The proof does not extend to
genus $2$, because some of the families used do not yield test
curves in the genus $2$ case but only points. (For example one
family is constructed by attaching a fixed elliptic curve to a
moving point on a fixed $g-1$ curve. For genus $g=2$ all the curves
in the family are isomorphic.).

\subsection{Hodge numbers}\label{ss9}

\begin{Thm} \label{4}
 For every $\overline{X} \in \{ \bR_2, \bS_2^+, \bS_2^- \}$, the rational cohomology of $\bR_2$ is algebraic,
i.e. all odd cohomology groups vanish,
and for all $n \in \mbb{N}$ we have $H^{2n}_{\mbb{Q}}( \overline{X}) = A^n_{\mbb{Q}} (\overline{X})$. Furthermore:

(i) The boundary classes generate the $\mbb{Q}$-vectorspace $H^2_{\mbb{Q}}( \overline{X})$.

(ii) There is an ample divisor $L$ which is a linear combination of
the boundary classes of $\overline{X}$, such that $L H^2_{\mbb{Q}}(
\overline{X}) = H^4_{\mbb{Q}}( \overline{X})$. Thus the products of
$L$ with the boundary classes generate the $\mbb{Q}$-vectorspace
$H^4_{\mbb{Q}}( \overline{X})$.

Hence the boundary classes generate the $\mbb{Q}$-algebras
$H^*_{\mbb{Q}}(\overline{X})$ and $A^*_{\mbb{Q}}(\overline{X})$.
\end{Thm}
{} \textbf{Proof:} All except part (ii) follows as a special case
from Corollary \ref{n3} (ii) and (iv).

Proof of (ii): $\bS_2$ being projective, there is an ample divisor
on this space. Like every divisor, according to lemma \ref{a4}, it
is equivalent to a linear combination $L$ of boundary classes. Of
course $L$ is also ample. According to the Hard Lefshetz Theorem,
multiplication with $L$ induces an isomorphism from $H^2_{\mbb{Q}}(
\overline{X})$ to $H^4_{\mbb{Q}}( \overline{X})$. The Hard Lefshetz
Theorem holds for our moduli spaces according to Summary \ref{p1} (i)
$\square$

\begin{Thm} \label{3}
 $\bR_2$, $\bS_2^+$ and $\bS_2^-$ all have Hodge diamonds of the following form

\[\begin{matrix} & & & 1
            \\ & & 0 & & 0
            \\ & 0 & & n & & 0
            \\ 0 & & 0 & & 0 & & 0
            \\ & 0 & & n & & 0
            \\ & & 0 & & 0
            \\ & & & 1

 \end{matrix}\]

with $n=4$ for $\bR_2$ and $n=3$ for $\bS_2^+$ as well as $\bS_2^-$.
\end{Thm}
{} \textbf{Proof:} For every $\overline{X} \in \{ \bR_2, \bS_2^+,
\bS_2^- \}$ $h^{2,0}(\overline{X}) = 0$ by Corollary \ref{n3} (v) ,
thus, due to the symmetries of the Hodge diamond, also
$h^{0,2}(\overline{X}) = 0$, $h^{1,3}(\overline{X}) = 0$ and
$h^{3,1}(\overline{X}) = 0$. Theorem \ref{4} then yields
$h^{1,1}(\overline{X})= h^{2,2}(\overline{X})$, and the value for
$n=h^{1,1}(\overline{X})$ is given by Lemma \ref{a4}. $\square$

\subsection{The cohomology rings in terms of generators and relations.}\label{ss10}

By Theorem \ref{4} we know that for our moduli spaces the Chow ring
and the rational cohomology ring coincide, and that they are
generated by the boundary classes. Now we determine the graded ring
structures:

\begin{Thm} \label{z}
 (i) The rational Chow ring $A^*_{\mbb{Q}}(\bR_2)$ is as a graded $\mbb{Q}$-Algebra isomorphic to the quotient
 $\mbb{Q}[\doi, \doii, \dor, \di, \dii]/I$, where $I$ is the homogeneous ideal generated by the following (independent) elements:

\[\doi + 6 \doii -3 \dor +12 \di -8 \dii, \]
\[\doii \dii, \qquad \doii \dor, \qquad \di \dii,\]
\[\di (\doii - \dor), \qquad \dii (\doi - \dor), \qquad 4 (\dii)^2 + \dor \dii,
 \qquad 2 \doi \doii + 4 \doi \di - 4 \doi \dii - \doi \dor , \]
\[\doi (\dor)^2, \qquad (\doi)^2 \doii \]

 (ii) $A^*_{\mbb{Q}}(\bS_2^+) \cong \mbb{Q}[\aop, \bop, \aip, \bip]/J$, where $J$ is the homogeneous ideal generated by the following (independent) elements:

\[3 \aop -4 \bop -8 \aip + 72 \bip,\]
\[\aip \bip, \qquad \bop \bip, \qquad \aop \aip - \bop \aip,  \]
\[(\aop)^2 \bop, \qquad (\aop)^2 (\aip - \bip) \]

 (iii) $A^*_{\mbb{Q}}(\bS_2^-) \cong \mbb{Q}[\aom, \bom, \aim]/K$, where $K$ is the homogeneous ideal generated by the following (independent) elements:

\[ 24 (\aim)^2 + \aom \aim + 2 \bom \aim, \qquad 12 (\bom)^2 +24 \bom \aim + \aom \bom , \]
\[3 (\aom)^2 -4 \aom\bom -8 \aom \aim + 80 \bom \aim\]

\end{Thm}
{}
\textbf{Proof:} The general idea of the proof and many of its steps are adopted from \cite{MR2526309}.

The rational Chow rings of our Moduli spaces are generated by the
boundary components according to Theorem \ref{4}. Thus there is a surjective morphism from the quotient algebras
of our Theorem to these Chow rings, if only the elements listed above as generators of the ideals of
relations $I$, $J$ and $K$, indeed equal zero in the rational Chow ring.

If this is shown, the following fact implies, that the morphisms are even isomorphisms:
The homogeneous components of the algebra $\mbb{Q}[\doi, \doii, \dor, \di, \dii]/I$
have $\mbb{Q}$-vectorspace dimensions $1,4,4,1,0,0,...$, whereas the
homogeneous components of $\mbb{Q}[\aop, \bop, \aip, \bip]/J$ and
$\mbb{Q}[\aom, \bom, \aim]/K$ have dimensions $1,3,3,1,0,0,...$, as one can check using a coputer algebra system like Macaulay 2.
These are exactly the vectorspace dimensions of the homogeneous components of the rational Chow rings (according to theorem \ref{3}).

To prove most of the relations, we will use the finite morphisms onto
boundary components described in section \ref{ss7}. By these
morphisms we will push forward classes and relations. Many of the relations we will push forward
are already described in section \ref{ss7}. Pushforwards of 
boundary cycles are listed in the tables of Lemma
\ref{a11}. In the computations we will use these facts without
mentioning that we take them from section \ref{ss7}.
\\
\\ First we prove the relations for $\bR_2$.

The linear relation

\begin{equation} \label{g0}
 \doi + 6 \doii -3 \dor +12 \di -8 \dii =0
\end{equation}

holds by Lemma \ref{a8}.

A Prym curve corresponding to a point in $\Doii$ can not
correspond to a point in $\Dii$. The preimage of such a point under
$\toii: \bM_{1,2} \longrightarrow \Doii$, would have to correspond
to a reducible curve. Such a curve is of the following  form: It
consist of a component $D$ of genus 1, and a component $E \cong
\mbb{P}^1$ with two marked points on it. $D$ and $E$ meet in one
node. The Prym curve generated by gluing the marked points has a
genus 1 component corresponding to $D$. Restricted to this component
its Prym sheaf is trivial. The Prym curve can thus not correspond to
a point in $\Dii$. So $\Doii \cap \Dii = 0$, and:

\begin{equation} \label{g2}
 \doii \dii = 0
\end{equation}

Similarly one can prove

\begin{equation} \label{g3}
 \doii \dor = 0
\end{equation}

and

\begin{equation} \label{g1}
 \di \dii = 0
\end{equation}

Now we use the morphism $\ti: \bM_{1,1} \times \bR_{1,1} \longrightarrow
\Di$. In $A_{\mbb{Q}}^1(\bR_{1,1})$ the relation $\tilde{\doii} =
\tilde{\dor}$ holds . Thus we also have $1 \otimes \tilde{\doii} = 1
\otimes \tilde{\dor}$ in $A_{\mbb{Q}}^1(\bM_{1,1} \times
\bR_{1,1})$. Pushing this forward by $\ti$ one gets:

\[(\ti)_*( 1 \otimes \tilde{\doii}) =  (\ti)_*(1 \otimes \tilde{\dor}) \]
\[ \Leftrightarrow \quad  \di \doii  = \di \dor \]

\begin{equation} \label{g4}
 \Leftrightarrow \quad \di(\doii - \dor) = 0
\end{equation}

Similarly, but using the  $\tii: \bR_{1,1} \times \bR_{1,1}
\longrightarrow \Dii$, we get:

\begin{equation} \label{g5}
 \dii(\doi - \dor) = 0
\end{equation}

According to \cite{MR717614}, page 321, in $A_{\mbb{Q}}^*(\bM_2)$
the relation $10 \lambda = \delta_0 + 2 \delta_1$ holds. Pulling
this back by $\pi_R$ to $\bR_2$ one gets:

\begin{equation} \label{g6}
 l = \frac{1}{10} (\doi + \doii + 2 \dor + 2 \di + 2 \dii)
\end{equation}

Multiplying equation (\ref{g6}) with $\dii$ and using equations (\ref{g1}), (\ref{g2}) and (\ref{g5}) yields:

\begin{equation} \label{g7}
 \dii l = \frac{1}{10} (3 \dii \dor + 2 (\dii)^2)
\end{equation}

On the other hand, because of $\dii = \frac{1}{2} (\tii)_*(1))$ we
can write $\dii l = \frac{1}{2} (\tii)_*((\tii)^* l)$ by the
projection formula. According to the Lemmata \ref{7} and \ref{8}

\[(\tii)^* l = \tilde{l} \otimes 1 + 1 \otimes \tilde{l} = \frac{1}{4} (\tilde{\dor} \otimes 1) + \frac{1}{4} (1 \otimes \tilde{\dor}) \]

We use $\dii \dor = (\tii)_* (\tilde{\dor} \otimes 1) =  (\tii)_* ( 1 \otimes \tilde{\dor})$ and get:
\[ \dii l = \frac{1}{2} (\tii)_*((\tii)^* l) = \frac{1}{2} (\tii)_* ( \frac{1}{4} (\tilde{\dor} \otimes 1) + \frac{1}{4} (1 \otimes \tilde{\dor})) \]
\[=  \frac{1}{2}  \frac{1}{4} (\dii \dor + \dii \dor) = \frac{1}{4} \dii \dor \]

By subtracting the equation $ \dii l = \frac{1}{4} \dii \dor$ from equation (\ref{g7}), and multiplying by 20, one gets:

\begin{equation} \label{g8}
 4 (\dii)^2 + \dor \dii = 0
\end{equation}

The last codimension $2$ relation

\begin{equation} \label{g9}
\qquad 2 \doi \doii + 4 \doi \di - 4 \doi \dii - \doi \dor
\end{equation}

we have proven earlier (Lemma \ref{a10}).

To obtain the codimension $3$ relations we use that $l^2 \doi = l^2 \doii = l^2 \dor = 0$ (cf. lemma \ref{a12}).

Because of $\doii= \frac{1}{2} (\toii)_*1$ we can write $\doii l = \frac{1}{2} (\toii)_*((\toii)^* l)$.  According to Lemma \ref{7} and \ref{8} one has

\[(\toii)^* l = \hat{\lambda} = \frac{1}{12} \hat{\delta_0}\]

By using $\doi \doii = \frac{1}{2} (\toii)_* \hat{\delta_0}$ we get

\[ \doii l = \frac{1}{2} (\toii)_*(\frac{1}{12} \hat{\delta_0}) = \frac{1}{12} \doi \doii\]

Thus $0= l^2 \doii = \frac{1}{12}l \doi \doii = \frac{1}{144} (\doi)^2 \doii$, and so

\begin{equation} \label{g11}
(\doi)^2 \doii = 0
\end{equation}

Using $\doi= \frac{1}{2} (\toi)_*1$ we can write $\doi l = \frac{1}{2} (\toi)_*((\toi)^* l)$.  According to Lemma \ref{7} and \ref{8} one has

\[(\toi)^* l = \bar{l} = \frac{1}{6} (c'+c''+c^r)\]

By using the pushforwards of Lemma \ref{a11} we get

\[ \doii l = \frac{1}{2} (\toi)_*(\frac{1}{6} (c'+c''+c^r)) = \frac{1}{12} (2 [E^{\prime,\prime}]_Q + \doi \doii + 2 \doi \dor)\]

Together with the relation $2 [E^{\prime,\prime}]_Q + \doi \doii =  \doi \dor$ of Lemma \ref{a10} (iii), this yields

\[ \doi l = \frac{1}{4} \doi \dor\]

Thus $0= l^2 \doi = \frac{1}{4}l \doi \dor = \frac{1}{16} \doi (\dor)^2$, and so

\begin{equation} \label{g10}
\doi (\dor)^2 = 0
\end{equation}

We have proven that the generators of the ideal $I$ are indeed equal to $0$ in the rational Chow ring of $\bR_2$.
\\
\\Now we prove the relations on $\bS_2^+$

The linear relation

\begin{equation} \label{h1}
 3 \aop -4 \bop -8 \aip + 72 \bip =0
\end{equation}

holds by Lemma \ref{a8}.

Similar to what was done for $\bR_2$ above, one can show that $\Aip
\cap \Bip = \emptyset$ and $\Bop \cap \Bip = \emptyset$, so we have
the relations

\begin{equation} \label{h2}
  \aip \bip = 0
\end{equation}

and

\begin{equation} \label{h3}
  \bop \bip = 0
\end{equation}

Proceeding similarly as in the proof of equation \ref{g4} and using
the morphism $\ria: \bS_{1,1}^+ \times \bS_{1,1}^+ \longrightarrow
\Aip$ we get:

\begin{equation} \label{h4}
 \aip(\aop - \bop) = 0
\end{equation}

To obtain the codimension $3$ relations, similar to the case of
$\bR_2$ we use that $\aop (l^+)^2  = \bop (l^+)^2  = 0$ (c.f. Lemma
\ref{a12}).

Because of $\bop= \frac{1}{2} (\rob)_*1$ we can write  $\bop l^+ =
\frac{1}{2}(\rob)_*((\rob)^* l^+)$. According to Lemma \ref{7} and
\ref{8} one has

\[(\rob)^* l^+ = \hat{l}^+ = \frac{1}{4} \hat{\alpha}_0^+\]

By using $\aop \bop = \frac{1}{2} (\rob)_* \hat{\alpha}_0^+ $ we get
\[ \bop l^+ = \frac{1}{2} (\toi)_*(\frac{1}{4} \hat{\alpha}_0^+) = \frac{1}{4} \aop \bop\]

Thus $0= \bop (l^+)^2 = \frac{1}{4} \aop \bop l^+ = \frac{1}{16} (\aop)^2 \bop$, and so

\begin{equation} \label{h5}
(\aop)^2 \bop = 0
\end{equation}

We would also like to make use of $\aop (l^+)^2=0$, by expressing
$\aop (l^+)^2$ in a nontrivial way as a product of boundary classes,
but the morphism $\roa$ does not help. We instead use equation (\ref
{h1}) to write  $3\aop$ as $4 \bop +8 \aip - 72 \bip$ and to get $0=
(4 \bop +8 \aip - 72 \bip)(l^+)^2$. Because of $\bop (l^+)^2=0$ this
simplifies to

\begin{equation} \label{h6}
(\aip -9 \bip) (l^+)^2 =0
\end{equation}

We can write $\aip l^+ = \frac{1}{4} (\ria)_*((\ria)^* l^+ )$, and here the Lemmata \ref{7} and \ref{8} yield

\[(\ria)^* l^+ = \tilde{l}^+ \otimes 1 + 1 \otimes \tilde{l}^+ =
\frac{1}{4}( \tilde{\alpha}_0^+  \otimes 1 + 1 \otimes \tilde{\alpha}_0^+)\]

By using
$\aop \aip = \frac{1}{2} (\ria)_* (\tilde{\alpha}_0^+  \otimes 1)=\frac{1}{2} (\ria)_* (1 \otimes \tilde{\alpha}_0^+)$
we get

\[ \aip l^+ = \frac{1}{4} (\ria)_*(\frac{1}{4} (\tilde{\alpha}_0^+  \otimes 1 + 1 \otimes \tilde{\alpha}_0^+)) = \frac{1}{4} \aop \aip \]

Analogously, from $\bip l^+ = \frac{1}{4} (\rib)_*((\rib)^* l^+ )$
we get to

\[ \bip l^+ = \frac{1}{4} (\rib)_*(\frac{1}{12} (\tilde{\alpha}_0^+  \otimes 1 + 1 \otimes \tilde{\alpha}_0^+)) = \frac{1}{12} \aop \bip \]

By using $\aip l^+ = \frac{1}{4} \aop \aip$ and $\bip l^+ =
\frac{1}{12} \aop \aip$ one can now rewrite equation (\ref{h6})

\[0 = (\aip -9 \bip) (l^+)^2 = \aop ( \frac{1}{4} \aip -9 \frac{1}{12}\bip) l^+ = (\aop)^2 ( \frac{1}{16} \aip -9 \frac{1}{144}\bip)\]

Thus

\begin{equation}
(\aop)^2(\aip - \bip) =0
\end{equation}

(The codimension $3$ relations computed in \cite{MR2526309},
except of $(\aop)^2 \bop = 0$, are incompatible with our results.)
\\
\\ Now we come to the relations on $\bS_2^-$.

The relation $12 (\delta_1)^2 + \delta_0 \delta_1 =0$ holds on
$\bM_2$ as follows directly from Theorem 10.1. of \cite{MR717614}.
Pulling this relation back by $\pi_-$ yields the first relation

\begin{equation} \label{i1}
24 (\aim)^2 + \aom \aim + 2 \bom \aim =0
\end{equation}

Pulling back the relation $10 \lambda = \delta_0 + 2 \delta_1$ by
$\pi_-$ one gets:

\begin{equation} \label{i2}
 l^- = \frac{1}{10} (\aom + 2 \bom + 4 \aim)
\end{equation}

Multiplication by $\bom$ yields:

\begin{equation} \label{i3}
 l^- \bom = \frac{1}{10} (\aom \bom + 2 (\bom)^2 + 4 \bom \aim)
\end{equation}

On the other hand, because of $\bom = \frac{1}{2} (\hob)_*(1))$, we can write $\bom l^- = \frac{1}{2} (\hob)_*((\hob)^* l)$.
According to the Lemmata \ref{7} and \ref{8}

\[(\hob)^* l^- = \hat{\lambda} = \frac{1}{12} \hat{\delta_0} \]

We use $\aom \bom = \frac{1}{2}(\hob)_* \hat{\delta_0}$ and get:
\[ l^- \bom = \frac{1}{2} (\hob)_*(\frac{1}{12} \hat{\delta_0}) = \frac{1}{12} \aom \bom \quad (\ast)\]

By subtracting the equation $ \bom l^- = \frac{1}{12} \aom \bom$ from equation (\ref{i3}), and multiplying by 60, one gets:

\begin{equation} \label{i4}
 12 (\bom)^2 +24 \bom \aim + \aom \bom
\end{equation}

(In \cite{MR2526309} it is claimed that $l ^- \bom =
\frac{1}{6} \aom \bom$ instead of $(\ast)$, from this then follows
$3 (\bom)^2 + 6 \bom \aim - \aom \bom$ instead of equation
(\ref{i4}).)

To get the last relation we first compute three relations containing
classes that can not immediately be written as products of boundary
classes (for the description of the closed strata defining these
classes, c.f. the appendix). The fist of these relations we take
from Lemma \ref{a10}:

\begin{equation} \label{i5}
16 [X^-]_Q + [C^-]_Q - 4 \aom \aim - \aom \bom = 0
\end{equation}

In $A_{\mbb{Q}}^1(\bS_{1,1}^+)$ the relation $\tilde{\alpha}_0^+ =
\tilde{\beta}_0^+$ holds, which implies for
$A_{\mbb{Q}}^1(\bS_{1,1}^+ \times \bM_{1,1})$ the relation
$\tilde{\alpha}_0^+ \otimes 1 = \tilde{\beta}_0^+ \otimes 1$.
Pushing this forward by the morphism $\hia: \bS_{1,1}^+ \times
\bM_{1,1} \longrightarrow \Aom \subset \bS_2^-$ yields:

\begin{equation}\label{i6}
[X^-]_Q = \bom \aim
\end{equation}

(In \cite{MR2526309} the authors claim, that one can get the
equation $\aom \aim = \bom \aim$ instead of equation (\ref{i6}).
Using the projection formula and the morphism $\hia$ they obtain the
equation $\aom \aim - (\hoa)_*(1 \otimes \delta_0) = \bom \aim$.
Then they claim that $(\hia)_*(1 \otimes \delta_0) = \frac{1}{2}
\aom \aim$, from which their equation would follow. If I understand
them correctly, they assume that $\bS_{1,1}^+ \times \Delta_0$ is
mapped $1:1$ onto $\Aom \cap \Aim$ by $\hia$. This would be wrong.
$\bS_{1,1}^+ \times \Delta_0$ is only mapped onto $Y^-$, which is one of the two
irreducible components of $\Aom \cap \Aim$, the other being $X^-$.
There is no a priori reason for $[Y^-]_Q$ and $[X^-]_Q$ to be
equivalent, so their equation does not follow. As one can check
after computing all relations, the equation does not hold.)

By multiplying equation (\ref{i2}) with $\aom$ one gets

\begin{equation} \label{i7}
 l^- \aom= \frac{1}{10} ((\aom)^2 + 2 \aom \bom + 4 \aom \aim)
\end{equation}

On the other hand, because of $\aom = \frac{1}{2} (\hoa)_*(1))$, we can write $\aom l^- = \frac{1}{2} (\hoa)_*((\hoa)^* l)$.
According to the Lemmata \ref{7} and \ref{8}

\[(\hoa)^* l^- = \check{l} = \frac{1}{12}(\check{\alpha}_0 + 2 \check{\beta}_0) \]

We use $[C^-]_Q = \frac{1}{4} (\hoa)_* \check{\alpha}_0$ and $\aom \bom = \frac{1}{2} (\hoa)_* \check{\beta}_0$ to get :
\[ l^- \aom = \frac{1}{2} (\hoa)_*(\frac{1}{12} (\check{\alpha}_0 + 2 \check{\beta}_0)) =
\frac{1}{6}([C^-]_Q + \aom \bom)\]

By subtracting the equation $ l^- \aom = \frac{1}{6}([C^-]_Q + \aom
\bom)$ from equation (\ref{i7}), and multiplying by 30, one gets:

\begin{equation} \label{i8}
 5 [C^-]_Q = 3 (\aom)^2 + \aom \bom + 12 \aom \aim
\end{equation}

Plugging equation (\ref{i6}) into equation (\ref{i5}) yields:

\[16 \bom \aim+ [C^-]_Q - 4 \aom \aim - \aom \bom = 0\]

By multiplying this with $5$ and plunging in equation (\ref{i8}) we get

\begin{equation} \label{i9}
3 (\aom)^2 -4 \aom\bom -8 \aom \aim + 80 \bom \aim
\end{equation}

This is the last relation we had to check. $\square$
\\
\\ \textbf{Remarks:} (i) One can test these relations by pulling the known relations $\delta_0 \delta_1 + 12 (\delta_0)^2=0$ and $528 (\delta_1)^3 + (\delta_0)^3 =0$
(known from \cite{MR717614}) back from $\bM_2$ to our moduli spaces
and check whether they are fulfilled in the rings that Theorem
\ref{z} claims to be to the rational Chow rings.

(ii) While the cohomology rings of $\bS_2^+$ and $\bS_2^-$ have,
according to our computation, the same Betti numbers, they are still
nonisomorphic as is clear by the fact that the relations in
codimension $1$ and $2$ determine the cohomology ring of $\bS_2^-$
completely, whereas for $\bS_2^+$, codimension $3$ relations are
needed.

\section{Appendix} \label{s5}

\subsection{Stratifications ``by topological type''}

We now describe the strata of the stratifications of $\bM_2$, $\bS_2^+$,
$\bS_2^-$ and $\bR_2$ according to the topological type of the
curves. For one of the moduli spaces beside $\bM_2$ we mean by this the
irreducible components of the preimages of the strata of $\bM_2$
under the forgetful morphism. In what follows we do rather describe
the closures of the strata than the strata themselves. We call these
closures the closed strata of the stratifications according to
topological type.

The description of the stratifications of $\bS_2^+$ and $\bS_2^-$
can be found in the appendix of \cite{MR2526309}. We use the symbols
introduced there for the strata, instead to denote the closures of
the strata, because some of the associated cycle classes appear in
our computations. The stratification of $\bM_2$ is described in
\cite{MR717614} §9., we will use the notation introduced there for
the closed strata.

We now describe the closed strata of $\bR_2$. We group the strata according to the
strata of $\bM_2$ they are lying over. For every stratum we
will explain how a Prym curve $(X;L;b)$ parametrized by a general
point looks like. We call such a Prym curve a general Prym curve of
the stratum.

Strata over $\Delta_0$: For these cycles the underlying stable model $C$ of a generic Prym curve is irreducible with one node.
The cycles are the divisors $\Doi$, $\Doii$ and $\Dor$ described in section \ref{ss2}

Strata over $\Delta_1$: For these cycles the underlying stable model
$C$ of a generic Prym curve consist of two smooth irreducible
components meeting in one node. Again these divisors ($\Di$ and
$\Dii$) are described in section \ref{ss2}

Strata over $\Delta_{00}$: For these codimension $2$ strata the
underlying stable model $C$ of a general Prym curve is an
irreducible curve with two nodes.
\begin{enumerate}
 \item $E^{\prime,\prime}$. General Prym curve: $X = C$, normalizing either of the two nodes and pulling back $\mc{L}$ to this partial normalization yields a nontrivial Prym sheaf.
 \item $E^{\prime,\prime \prime}$. General Prym curve: $X = C$, normalizing one of the two nodes and pulling back $\mc{L}$ to this partial normalization yields in one case a nontrivial Prym sheaf, in the other case the trivial sheaf, depending on which node was normalized .
 \item $E^{\prime,r}$. General Prym curve: $X$ is obtained from $C$ by blowing up one of the two nodes.
 \item $E^{r,r}$. General Prym curve: $X$ is obtained from $C$ by blowing up both nodes.
\end{enumerate}

Strata over $\Delta_{01}$: For these codimension $2$ strata the
underlying stable model $C$ of a general Prym curve consists of two
irreducible components, one of them, called $C_1$, is smooth, the
other one, $C_2$,  has a node.
\begin{enumerate}
  \item $F_1'$. General Prym curve:  $X=C$, $\mc{L}_{|C_1}$ is nontrivial, $\mc{L}_{|C_2}$ is trivial.
  \item $F_1''$. General Prym curve:  $X=C$, $\mc{L}_{|C_1}$ is trivial, $\mc{L}_{|C_2}$ is nontrivial.
  \item $F_1^r$. General Prym curve:  $X$ is obtained from $C$ by blowing up the node on $C_2$, $\mc{L}_{|C_1}$ is trivial.
  \item $F_{1:1}'$. General Prym curve: $X=C$, both $\mc{L}_{|C_1}$ and $\mc{L}_{|C_2}$ are nontrivial.
  \item $F_{1:1}^r$. General Prym curve: $X$ is obtained from $C$ by blowing up the node on $C_2$, $\mc{L}_{|C_1}$ is nontrivial.
\end{enumerate}

Strata over $C_{000}$: For these codimension $3$ strata the underlying stable model $C$ of a general Prym curve
consists of two irreducible smooth rational components meeting in three nodes.
\begin{enumerate}
  \item $G'$. General Prym curve:  $X=C$
  \item $G^r$. General Prym curve: $X$ is obtained from $C$ by blowing up one of the nodes.
\end{enumerate}

Strata over $C_{001}$: For these codimension $3$ strata the underlying stable model $C$ of a general Prym curve
consists of two irreducible components $C_1$ and $C_2$ meeting in one node, each irreducible component having one node.
\begin{enumerate}
  \item $H_1'$. General Prym curve: $X=C$, restricting $\mc{L}$ to one of the components yields a nontrivial Prym sheaf, restricting to the other yields the trivial sheaf.
  \item $H_1^r$. General Prym curve: $X$ is obtained from $C$ by blowing up the node on one of the components, $\mc{L}$ is trivial restricted to the component not blown up.
  \item $H_{1:1}'$. General Prym curve: $X=C$, $\mc{L}$ is nontrivial on both components.
  \item $H_{1:1}^r$. General Prym curve: $X$ is obtained from $C$ by blowing up the node on one of the components, $\mc{L}$ is nontrivial on both components.
  \item $H_{1:1}^{r,r}$. General Prym curve: $X$ is obtained from $C$ by blowing up the nodes on both components.
\end{enumerate}

\subsection{Comparison of automorphisms}

As we have shown, there is an isomorphism of coarse moduli spaces
$a_R:\bM_{0,[2,4]} \overset{\cong}{\rightarrow} \bR_2$, and also
$\bS_2^+$ and $\bS_2^-$ are isomorphic to moduli spaces of stable
genus $0$ curves with partitioned marked Points. Now, lets say for
$x \in \bM_{0,[2,4]}$, one can ask how the Automorphisms of
objects in the class $x$ and objects in the class $a_R(x)$ fit
together. 

The fact that we know $a_R$ explicitly only on an open subset of
$\bM_{0,[2,4]}$ makes it difficult to compare the automorphisms on
both sides directly. But one can overcome this difficulty by
extending the automorphisms to the local universal deformation
spaces of the Prym curve belonging to $a_R(x)$, respectively to the local universal deformation space 
of the stable genus $0$ curve with marked points belonging to $x$. This is
helpful because on the loci of smooth curves of these
deformation spaces, i.e. almost everywhere, one knows explicitly how Prym
structures and marked points correspond to each other. 

(We only spoke about $a_R$ in this introduction, but
everything also applies to $a_+$ and $a_-$.)

\textbf{General notation (part 1):} For the rest of the section let
$\bM_{0,\bullet}$ denote one of the moduli spaces $\bM_{0,[1,5]}$,
$\bM_{0,[2,4]}$ and $\bM_{0,[3,3]}$. We denote by $\bQ_2$ the one of
the moduli spaces $\bR_2$, $\bS_2^+$ and $\bS_2^-$ the space
$\bM_{0,\bullet}$ is isomorphic to. We call $a_{\bullet}:
\bM_{0,\bullet} \rightarrow \bQ_2$ the corresponding isomorphism
(one out of $a_R$, $a_+$ and $a_-$). (C.f. Lemma \ref{2})

If not specified otherwise $(D;\{A,B\})$ is always a genus $0$ curve
$D$ together with two disjoint sets $A,B$ of marked points, such that
$(D;\{A,B\})$ is parametrized by a point $x \in \bM_{0,\bullet}$ . We
define $y:= a_{\bullet}(x)$.

$f: Y \rightarrow D$ is always the admissible $2:1$ cover of $(D; A
\cup B)$.

\begin{Def}
 If $(D;M)$ is a stable genus $0$ curve with a set $M$ of marked points, then we call those irreducible components of $D$
the \emph{extremities of} $(D;M)$ which meet the rest of $D$ only in
one point and which carry only two of the marked points.
\end{Def}

\begin{Lemma} \label{ap4}

(i) For $f: Y \rightarrow D$ as in the general assumption, $E$ an
extremity of $(D;A \cup B)$, the preimage $f^{-1}(E)$ is an
exceptional component of the quasistable curve $Y$.

(ii)  Let $Cont_1: Y \rightarrow C$ be the contraction of all
extremal components of $Y$, $C$ the stable model of $Y$. Let
$cont_1: D \rightarrow \hat D$ be the contraction of all extremities
of $(D; A \cup B)$. Then there is a unique finite $2:1$ morphisms $C
\rightarrow \hat D$ fitting into the following commutative diagram.
\[ \begin{xy} \xymatrix{ Y \ar[r]^{Cont_1}  \ar[d]_f & C \ar[d]_{\hat f} \\
                         D \ar[r]_{cont_1}          & \hat D } \end{xy}\]

(iii) There is a (not necessarily unique) way to blow up nodes of
the curve $C$ such that the variety $X$ obtained by this can be
equipped with a structure $(\mc{L};b)$ such that $[(X; \mc{L}; b)] =
y$. We pick such a blowup morphism $Cont_2: X \rightarrow C$.
\end{Lemma}

\textbf{Proof:} (i): Look at the explicit description of an
admissible $2:1$ cover in the proof of Cor. 2.5 in \cite{MR1981190}.

(ii): For every exceptional component $E$ contracted to a point by
$Cont_1$, $cont_1$ contracts the extremity $f(E)$ to a point.
Considering this, it is obvious how to define $\hat f: Y \rightarrow
\hat D$ to fit the diagram.

(iii): By the construction of $a_{\bullet}$ in the proof of Lemma \ref{n1}, it is clear that
\[ \begin{xy} \xymatrix{ \bM_{0,\bullet} \ar[r]^{a_{\bullet}}  \ar[d]_{\pi} & \bQ_2 \ar[d]_{\pi'} \\
                         \bM_{0,[6]} \ar[r]_{b}          & \bM_2 } \end{xy}\]
commutes, and thus for any $(X; \mc{L};b)$ with $[(X;\mc{L};b)]=y$, $X$ has the same
stable model $C$ as $Y$. $\square$

\textbf{General notation (part 2):} We will keep the notation of
Lemma \ref{ap4} for the rest of the section. We fix one $Cont_2: X
\rightarrow C$ as in part (iii).

We will denote by $\hat A$ resp. $\hat B$ the set of all points of
$\hat D$ that come from those marked points in $A$ resp. $B$ that
lie on components of $D$ not contracted by $cont_1$. By $H$ we will
denote the set of points of $\hat D$ to which extremities of $(D;A
\cup B)$ are contracted by $cont_1$. We set $G:= \hat A \cup \hat
B$. The object $(\hat D; (G,H))$ is then stable and parametrized by
a point of some $\bM_{0,(i_1,i_2)}$ with $3 \le i_1+i_2 \le 6$. If
we want to retain more information about the extremities contracted
to points of $H$, we decompose this set into $H_A$, $H_B$, $H_{AB}$,
where $H_A$ contains the points to which extremities carrying only
marked points of $A$ are contracted, $H_B$ contains those coming
from extremities with marked point only from $B$, while to the
points of $H_{A,B}$ extremities that carry one point of $A$ and one
point of $B$ are contracted. Note that then $( \hat D; \{(\hat A,
H_A), (\hat B, H_B)\},H_{AB})$ contains the full information about
the isomorphism class of $(D; \{A,B\})$.

\begin{Lemma} \label{ap3}

Using the general notation for this section:

(i) For every $\varphi \in Aut(Y)$ and every $a,b \in \tilde Y$
($\tilde Y$ the non-exceptional subcurve of $Y$):
\[f (a) = f (b) \Leftrightarrow  f (\varphi (a)) = f (\varphi(b))\]
Thus for every $\varphi \in Aut(C)$ and every $a,b \in C$:
\[\hat f (a) = \hat f (b) \Leftrightarrow  \hat f (\varphi (a)) = \hat f (\varphi(b))\]

(ii) There are natural surjective group homomorphisms
\[ Aut(Y) \overset{\chi_1}{\longrightarrow} Aut (C) \overset{\chi_2}{\longrightarrow} Aut((\hat D; (G,H))\]
and
\[Aut(X) \overset{\chi_1'}{\longrightarrow} Aut(C)\]
There is also a natural surjective group homomorphism
\[\psi_1': Aut( (D; A \cup B)) \rightarrow Aut((\hat D; (G,H)))\] (The homomorphisms
are defined explicitly in the proof.)

(iii) The kernels $Ker \, \chi_1$ resp. $Ker \, \chi_1'$ consist of
those automorphisms that are nontrivial only on the exceptional
components of $Y$ resp. $X$. The kernel $Ker \, \psi_1'$ consists of
those $\varphi \in Aut((D;\{A,B\})$ that are nontrivial only on the
extremities of $(D; A \cup B)$.

(iv) Let $\hat D_i$ be an irreducible component of $\hat D$, $C_i:=
\hat f^{-1}(\hat D_i)$, and let $\hat f_i: C_i \rightarrow \hat D_i$
be the restriction of $\hat f$. There is an automorphism $h_i$ on
$C_i$ interchanging the two sheets of $\hat f_i$. We call $h_i$ the
hyperelliptic involution on $C_i$ (although $C_i$ may be reducible
and not a hyperelliptic curve in the usual sense). One can extend
$h_i$ to an Automorphism of $C$, such that it is the identity on all
components of $C$ except $C_i$. We again denote this extension by
$h_i$.

(v) The $h_i \in Aut(C)$ belonging to the different irreducible
components of $\hat D$ generate the kernel $Ker \, \chi_2$. We call
the unique automorphism of $Ker \, \chi_2$ whose restriction to no
component of $C$ is trivial the \emph{full} hyperelliptic involution
of $C$
\end{Lemma}

\textbf{Proof:} (i): As shown in \cite{MR1981190}, the $2:1$
admissible cover of a stable genus $0$ curve $D$ with $2g+2$ marked
points for $g \ge 2$ is unique up to isomorphism. There also an
explicit method is given to associate to such a $D$ an admissible
$2:1$ cover. We use this explicit description in our proof, and one might
need to know it in order to understand the arguments.

Let $\tilde D$ be the subcurve of $D$ consisting of all components
of $D$ that are no extremities of $(D;A \cup B)$, and let $\tilde D_1,..., \tilde D_m$
be the irreducible components of $\tilde D$. For $i=1,...m$, let
$D_i$ be the subcurve of $D$ consisting of $\tilde D_i$ and the
extremities of $D$ attached to $\tilde D$, and let $Y_i$ be the part
of $Y$ lying over $D_i$. Define $\tilde Y_i:= \tilde Y \cap Y_i$
, and denote by $\tilde f_i: \tilde Y_i \rightarrow D_i$ the
restriction of $f$ to $\tilde Y_i$.

Let each of $q_{i,1},...,q_{i,l_i}$ be a point in which $D_i$ meets
one other component $D_j$ of $D$, let $Q_{i,1},...,Q_{i,l_i}$ be the
sets of points in $Y$ lying over $q_{i,1},...,q_{i,l_i}$. Each $Q_{i,j}$ contains
one or two points, an is contained in $\tilde Y_i$.
Let $p_{i,1},...,p_{i,k_i}$ be the points on $D$ in which $\tilde{D}_i$ meets extremities, 
and let $P_{i,1},...,P_{i,k_i}$ be the sets of points in $Y$ lying
over $p_{i,1},...,p_{i,k_i}$. For any $\varphi \in Aut(Y)$, $f \circ
\varphi: Y \rightarrow D$ is again an admissible $2:1$ cover of $(D;
\{A,B\})$. For any $\tilde Y_i$, $\varphi_{|\tilde Y_i}$ is an
isomorphism of $(\tilde Y_i; \{Q_{i,1},..., Q_{i,l_i}\},
\{P_{i,1},..., P_{i,k_i}\})$ to some $(\tilde Y_j; \{Q_{j,1},...,
Q_{j,l_j}\}, \{P_{j,1},..., P_{j,k_j}\})$, where $j=i$ is possible.

Now we prove (i) by checking several cases separately.

1. If $\tilde Y_i$ is a smooth connected curve of genus $\ge 2$, the
assertion of (i) (and (ii)) holds, since then $\tilde Y_i$ has to be
hyperelliptic and every hyperelliptic curve has a \emph{unique}
$g_2^1$ (C.f. \cite{MR0463157}, Chapt. IV, Prop. 5.3.)

2. If $\tilde Y_i$ is a smooth connected curve of genus $1$, then,
for stability reasons, either $l_i >0$ or $k_i >0$. But knowing one
fiber of the induced $2:1$ cover, determines one $g_2^1$ on an
elliptic curve uniquely. (C.f. \cite{MR0463157}, Chapt. IV, § 4.)

3. If $\tilde Y_i$ is a smooth connected curve of genus $0$, then
$\tilde f_i: \tilde Y_i \rightarrow \tilde D_i$ is ramified in
exactly two points. Thus $\tilde D_i$ carries at most $2$ points of
$A \cup B$. Thus, for reasons of stability and because $\tilde D_i$
is not an extremity of $D$, $l_i + k_i \ge 2$ has to hold, and one
of the $P$'s and $Q$'s on $\tilde Y_i$ has to contain two elements.
So we know that (i) holds for two fibers of $\tilde f_i: \tilde Y_i
\rightarrow \tilde D_i$. Because knowing the behavior of an
isomorphism of $\mbb{P}^1$'s in $3$ points determines the
isomorphism, one can quite easily conclude from this that (i)
holds for all of $\tilde Y_i$.

4. Otherwise $\tilde Y_i$ consists of two connected components,
which both are smooth genus $0$ curves. In this case there are no
points of $A \cup B$ lying on $\tilde D_i$, for otherwise $\tilde
f_i :\tilde Y_i \rightarrow \tilde D_i$ was ramified there. Thus,
for $\tilde D_i$ to be stable, we must have $l_i + k_i \ge 3$. So we
know three fibers of $\tilde f_i$ for which (i) holds. Again it
easily follows that (i) holds for all of $\tilde Y_i$.

(ii): $\chi_1$ is defined by $\varphi \mapsto \varphi^*$ for every
$\varphi \in Aut(Y)$ where $\varphi^* \in Aut(X)$ is the
automorphism defined by $\varphi^*(x):=
Cont_1(\varphi((Cont_1^{-1}(x))))$ for all $x \in X$. $\chi_2$ is
defined analogously. They are surjective because every automorphism
of a curve can obviously be extended to any curve obtained from it by blowing up nodes. 
We define $\chi_3$ by $\varphi \mapsto
\varphi^*$ for every $\varphi \in Aut(C)$ where $\varphi^* \in
Aut(\hat D)$ is the automorphism defined by $\varphi^*(x):= \hat
f(\varphi((\hat f^{-1}(x))))$. The definition of $\varphi^*(x)$
indeed gives a point by (i). We have to check $\varphi^* \in
Aut((\hat D;(G,H))$: $\varphi^*$ maps points in $H$ to points in
$H$, because they correspond to the $p$'s introduced in the proof of
(i), and we saw there that $\varphi$ maps the points lying over them
again to such points. The points of $C$ lying over $G$ are exactly
the smooth ramification points of $\hat f$, and by (i) $\varphi$
has to map such points to such.

That $\chi_2$ is surjective can be proven by again using the
decomposition of $Y$ used in the proof of (i) and again checking it in
the four possible cases distinguished there.

The morphism $\psi_1$ obviously exists and is surjective.

(iii): Clear from the definition of $\chi_1$ and $\chi_2$.

(iv): Obviously $h_i$ exists (uniquely). $C$ is of genus $2$ and the
components $C_i$ are of genus $1$ at least, so there can be only two
of them, and they can meet only in one point. Thus one can extend
each $h_i$ to the other component by the identity.

(v): The Kernel of $\chi_3$ consists of all $\varphi \in Aut(C)$
such that $\hat f(\varphi(a))= \hat f(a)$ for all $a \in C$. Quite
obviously the $h_i$ generate this group. $\square$



\begin{Lemma} \label{ap1}
Let $x \in \bM_{0,\bullet}$ be a point parametrizing $(D; \{A,B\})$ and $y := a_{\bullet} (y)$ its image in $\bQ_2$.
Let $(X; \mc{L}; b)$ be a object parametrized by $y$. Then:

(i) $Aut((D;\{A,B\}))$ is a subgroup of $Aut((D; A \cup B))$ and we call the restriction of the Morphism  $\psi_1'$ of Lemma
\ref{ap3}
\[ \psi_1: Aut((D;\{A,B\})) \rightarrow Aut((\hat D;(G,H))) \]

$Aut((X; \mc{L}; b))$ is a subgroup of $Aut(X)$. We call the restriction of $\chi_2 \circ \chi_1'$ to this subgroup
\[ \psi_2: Aut((X; \mc{L}; b)) \rightarrow Aut((\hat D; (G,H))) \]

From now on we use the abbreviations $M:= Aut ((D; \{A,B\}))$ and $N:= Aut((X; \mc{L}; b))$.

(ii) The group $Aut (( \hat D; \{(\hat A, H_A), (\hat B,
H_B)\},H_{AB}))$, (for the definition of this, c.f. the general
notation (part 2) for this section), is a subgroup of $Aut (( \hat
D; G,H))$ and:
\[\psi_2(N) = \psi_1(M) = Aut (( \hat D; \{(\hat A, H_A), (\hat B, H_B)\},H_{AB})) \]

(iii) Let $r$ be the number of extremities of $(D; A \cup B)$, let
$r'$ be the number of those extremities whose two marked points
either lie both in $A$ or lie both in $B$. Let $s$ be the number of
irreducible components of $D$. Let $I$ be the group of irrelevant
automorphisms of $(X; \mc{L}; b)$. We define $h:= \# \psi_1 (M)= \#
\psi_2 (N)$, $i:= \# I$, $m:= \# M$ and $n:= \#N$, then:
\[ m= 2^{r'} \cdot h, \qquad n= 2^{(s-r)} \cdot i \cdot h\]
and thus
\[n = 2^{(s-r-r')} \cdot i \cdot m\]

One can also write $i$ as $2^{u-1}$ where $u$ is the number of
connected components of $\tilde{X}$ the non-exceptional subcurve of
$X$.
\end{Lemma}

\textbf{Proof:} The different assertions that one automorphism group
is a subgroup of another one, made in parts (i) and (ii), are all
quite obvious.

The first thing we prove is the first equation of part (ii).

We are in the situation described by the following commutative diagram.
 \[ \begin{xy} \xymatrix{ & Y \ar[d]^f  \ar[rd]^{Cont_1} &  & \ar[ld]^{Cont_2}  X & \ar@{~>}[l] \mc{L}, \, b  \\
                          \{A,B\} \ar@{~>}[r]& D \ar[rd]_{cont_1} &  C \ar[d]^{\hat f} &  & \\
                          &(G,H) \ar@{~>}[r]  & \hat D &   &} \end{xy}\]
Where the curly arrows are meant to symbolize that 
additional structures are attached to $D$, $\hat D$ and $X$

From now on, we will work in the category of complex analytic spaces.

Let $(\mc{D} \rightarrow S; \{\mc{A}, \mc{B}\})$ be the local
universal deformation of $(D; \{A,B\})$. Part of the deformation is
a identification $\psi$ of the ``central fiber'' lying over a
special point $s_0 \in S$ with $(D; \{A,B\})$. $\mc{A}$ and $\mc{B}$
are then sets of sections on $\mc{D}$ meeting the central fibers in
the points of $A$ resp. $B$. Possibly after making a base change on
$S$ we can extend the $2:1$ admissible cover $f: Y \rightarrow D$ to
a local universal deformation $\mathbf{f}: \mc{Y} \rightarrow
\mc{D}$ over $S$. I.e. $\mathbf{f}$ restricted to the fiber over
$s_0$ can be identified with $f$, in a way compatible with $\psi$.
To see that this is possible, c.f. \cite{MR664324}, Page 61-62.
There the same thing is done for the local universal deformation of
a stable genus $0$ curve with \emph{ordered} marked points, but, since the
elements of $\mc{A}$ and $\mc{B}$ are sections (not multi-sections)
on our local deformation space, we can just put an arbitrary
ordering on them and by this make $(\mc{D} \rightarrow S; \{\mc{A}, \mc{B}\})$ into a local universal
deformation of a curve with ordered marked points.

If we denote the morphism contracting the exceptional components of $\mc{Y}$
by $\mathbf{Cont_1}: \mc{X} \rightarrow \mc{C}$ and the one contracting the
extremities of $(\mc{D}, \mc{A} \cup \mc{B})$ by $\mathbf{cont_1}: \mc{D}
\rightarrow \hat{\mc{D}}$, there is, analogously to Lemma \ref{ap4} (ii),
a morphism of families $\mathbf{\hat f}: \mc{C} \rightarrow
\tilde{\mc{D}}$ forming a commutative diagram with
$\mathbf{Cont_1}$, $\mathbf{cont_1}$ and $\mathbf{f}$.

Blowing up the appropriate nodes of $C$ (c.f. Lemma \ref{ap4}
(iii))) and the loci in the deformation $\mc{C}$ to which these
nodes extend, we arrive at an isomorphism $ \mathbf{Cont_2}: \mc{X}
\rightarrow \mc{C}$, such that $\mc{X} \rightarrow S$ is a
deformation of $X$. After making a base change we can extend
$(\mc{L};b)$ to a Prym- resp. spin structure $(\mathbf{L};
\mathbf{b})$ on $\mc{X}$. (C.f. \cite{MR1082361}, Page 570.)

Now we have deformations over $(S,s_0)$, forming the diagram
\[ \begin{xy} \xymatrix{ & \mc{Y} \ar[d]^{\mathbf{f}}  \ar[rd]^{\mathbf{Cont_1}} &  & \ar[ld]^{\mathbf{Cont_2}}  \mc{X} & \ar@{~>}[l] \mathbf{L}, \, \mathbf{b}  \\
                          \{\mc{A},\mc{B}\} \ar@{~>}[r]& \mc{D} \ar[rd]_{\mathbf{cont_1}} &  \mc{C} \ar[d]^{\hat{ \mathbf{f}}} &  & \\
                          &  & \hat{\mc{D}} &   &} \end{xy}\]
And by restricting these families to the fibers over $s_0$ we get
back to the diagram above.


We now prove the fist equation of part (ii) of our Lemma. We have to show that for $\varphi \in Aut((\hat D;(G,H)))$,

\[\varphi \in Im (\psi_1) \Leftrightarrow \varphi \in Im (\psi_2) \qquad (\ast) \]

In what follows, we use Lemma \ref{ap3} (ii). First lift an
automorphism $\varphi \in Aut((\tilde D;(G,H)))$ to an element of
$Aut(f:Y \rightarrow D)$, where we denote by $Aut({f}: Y \rightarrow
D)$ the automorphisms of the admissible $2:1$ cover, i.e.
automorphisms $\psi$ of $Y$ for which there exists an automorphism
$\psi$ of $(D; A \cup B)$ such that $\psi' \circ f = f \circ \psi$.
We extend this automorphism to an automorphism $\check{\varphi} \in
Aut(\mathbf{f}: \mc{Y} \rightarrow \mc{D})$, which is possible
because $\mathbf{f}: \mc{Y} \rightarrow \mc{D}$ is a local universal
deformation of the $2:1$ admissible cover $f:Y \rightarrow D$. By
the definition of $Aut(\mathbf{f}: \mc{Y} \rightarrow \mc{D})$,
$\check{\varphi}$ induces an automorphism of $(\mc{D}; \mc{A} \cup
\mc{B})$ we call $\check{\varphi_1}$. We restrict $\check{\varphi}$
to $\mc{C}$ and lift it to an automorphism $\check{\varphi}_2 \in
Aut(\mc{X})$. Restricting $\check{\varphi}_1$ and $\check{\varphi}_2$
to the central fibers, yields automorphism $\varphi_1 \in Aut((D; A \cup
B))$ and $\varphi_2 \in Aut(X)$ which are liftings of $\varphi \in
Aut((D;(G,H)))$ via the homomorphisms of Lemma \ref{ap3} (ii).

We have $\varphi \in Im(\psi_1)$ iff $\varphi_1$ respects the
structure $\{A, B\}$, and $\varphi \in Im(\psi_2)$ iff $\varphi_2$
respects the structure $(\mc{L}; b)$. Respecting the structure here
means that $\{\varphi_1(A), \varphi_1(B)\} =\{A, B\}$, respectively
$\varphi_2^* \mc{L} \cong \mc{L}$, compatible with $b$. This in turn is equivalent to
$\check{\varphi}_1$ respecting $\{\mc{A}, \mc{B}\}$, respectively
$\check{\varphi}_2$ respecting $(\mathbf{L}; \mathbf{b})$.

Let $S'$ be the open dense subset of the base space $S$ over which all
fibers of $\mc{D} \rightarrow S$ are smooth. Then over $S'$ the
contraction morphisms $\mathbf{Cont_1}$, $\mathbf{Cont_2}$ and
$\mathbf{cont_1}$ are all isomorphisms, and the diagram of
deformations above collapses to :
\[ \begin{xy} \xymatrix{  &   \mc{C}' \ar[d]^{\mathbf{f}'=\hat{\mathbf{f}}'} & \ar@{~>}[l] \mathbf{L}', \, \mathbf{b}'  \\
  \{\mc{A}',\mc{B}'\} \ar@{~>}[r]& \mc{D}'  &  } \end{xy} \]
(the $'$ indicating restriction to the preimages of $S'$)

On $\mc{C}'$ we can also define a spin- resp. Prym structure in the
following way: Define a divisor $E$ as such a linear combination of the
preimages under $\mathbf{f}'$ of the sections in $\mc{A}'$ as
is described (for points) in Lemma \ref{a6} (i)-(iii), and then let
$\mathbf{L}''$ be the line bundle on $\mc{C}'$ coresponding to $E$.
One can show that $\mathbf{L}' \cong \mathbf{L}''$ using that restricted to the fiber over any point of $S'$, $\mathbf{L}''$ and $\mathbf{L}'$ are isomorphic, and using that $(\mc{X}; \mathbf{L};
\mathbf{b})$ is the local universal deformation of every one of its
fibers (c.f. \cite{MR1082361} Page 574).

$\check{\varphi}_1$ and $\check{\varphi}_2$ can also be restricted
to the preimages of $S'$, where we denote them by
$\check{\varphi}_1'$ and $\check{\varphi}_2'$. Now
$\check{\varphi_1}$ respects $\{\mc{A}, \mc{B}\}$ iff
$\check{\varphi}_1'$ respects $\{\mc{A}', \mc{B}'\}$, while
$\check{\varphi}_2$ respects $(\mathbf{L}; \mathbf{b})$ iff
$\check{\varphi}_2'$ respects $(\mathbf{L}'; \mathbf{b}')$.

But $\mathbf{L}'$ as shown above is just the line bundle coresponding
to $E$, and looking at the definition of the divisor $E$ and at Lemma \ref{a6} one sees 
that its class does not change under permutations of the sections in $\mc{A}' \cup \mc{B}'$ wich respect the
partition into $\mc{A}'$ and $\mc{B}'$.
So $\check{\varphi}_2'$ respects $(\mathbf{L}'; \mathbf{b}')$ iff
$\check{\varphi}_1'$ respects $\{\mc{A}', \mc{B}'\}$.

Going back in our chain of equivalences of ``respecting'' conditions,
the previous sentence translates to, $\varphi_2 \in Aut((X;\mc{L};b)
\Leftrightarrow \varphi_1 \in Aut((D; \{A, B \}$, which  implies the
equivalence $(\ast)$ we wanted to prove.

That also the second equation of Part (ii) holds is easy to check,
considering how $(\tilde D; \{(\tilde A, H_A), (\tilde B,
H_B)\},H_{AB})$ contains the information which kinds of marked points the
contracted extremities of $D$ carried.

(iii): We know $Ker \, \psi_1 = Ker \, \psi_1' \cap Aut((D;
\{A,B\}))$. By Lemma \ref{ap3} (iii) this means that $Ker \, \psi_1$
consists of those Automorphisms of $(D;\{A,B\})$ that are nontrivial
only on Extremities of $(D; A \cup B)$. For every such extremity
carrying marked points only from the set $A$ or only from the set
$B$, there is an Automorphism of $(D; \{A,B\})$ that swaps the two
marked points and is trivial away form the extremity. These
automorphisms generate $Ker \, \psi_1$ which consist thus of
$2^{r'}$ elements. This together with (ii) implies:
\[h = m/2^{r'} \quad \Leftrightarrow \quad m = 2^{r'} \cdot h \]
To get the next equation we use
\[\# (Ker \, \psi_2) = \# \Bigl ( Ker \, \chi_1' \cap Aut((X; \mc{L}; b)) \Bigr ) \cdot
\# \Bigl (Ker \, \chi_2 \cap \chi_1' \bigl( Aut((X; \mc{L}; b))
\bigr ) \Bigr) \] Considering Lemma \ref{ap3} (iii) and the
definition of the irrelevant automorphisms of $(X; \mc{L}; b)$ (c.f.
preliminaries), we see that $Ker \, \chi_1' \cap Aut((X; \mc{L};
b))$ is just the group of irrelevant automorphisms. Since the
``hyperelliptic involutions'' generating $Ker \, \chi_2$ (c.f. Lemma
\ref{ap3} (v)) act trivially on all Prym- or spin sheaves, $Ker \,
\chi_2$ is contained in $\chi_1' \bigl( Aut((X; \mc{L}; b)) \bigr
)$. By Lemma \ref{ap3} (iv), $\#(Ker \, \chi_2)= 2^{s-r}$. This
implies:
\[h = n/2^{s-r} \quad \Leftrightarrow \quad n = 2^{s-r} \cdot i \cdot h \]

For the last assertion of (iii), c.f. \cite{MR2551759} Prop. 2.7., in the case of spin curves. (There the number of irrelevant automorphisms is $2^{u}$ instead of $2^{u-1}$ due to the different definition of
automorphisms). For Prym curves c.f. \cite{MR2639318} Remark 6.3.
$\square$

\subsection{Automorphism numbers}

\begin{Lemma} \label{ap2}
Let $p_1,...,p_n$ be $n$ distinct points of $\mbb{P}^1$ in general
position. We describe, for different $n \in \mbb{N}$, the group
$A:=Aut(\mbb{P}^1;\{p_1,...,p_n\})$ of automorphisms of $\mbb{P}^1$
that map points of the set $\{p_1,...,p_n\}$ again to points of this
set.

(i) For $n \le 2$, $A$ is an infinite group.

(ii) For $n=3$, $A$ has $6$ elements corresponding to the permutations of the $3$ points.

(iii) For $n=4$, $A$ has $4$ elements, one is the identity, the
others correspond to choosing two disjoint pairs of the points,
and interchanging the points in each pair.

(iv) For $n \ge 5$, $A$ consists only of the identity.
\end{Lemma}

\textbf{Proof:} The automorphisms of $\mbb{P}^1$ are the
transformations $x \mapsto \frac{Ax+B}{Cx+D}$ for $A,B,C,D \in
\mbb{C}$. With this information one can check that the assertions of
the Lemma are true. $\square$.

We can use Lemma \ref{ap2} together with Lemma \ref{ap1} (iii) to
compute the number of automorphisms of a general Prym- or spin curve
of one of the strata of the stratifications by topological type. Lemma \ref{ap1} allows to reduce to
computing the automorphism number of the corresponding genus $0$
curve with $6$ sorted marked points. The diagrams of these
corresponding objects are listed in the table below, and using Lemma
\ref{ap2} their automorphisms can quite easily be counted.

It may be even easier to draw a diagram for the object $(\tilde D;
\{(\tilde A, H_A), (\tilde B, H_B)\},H_{AB})$, count the number $h$
of automorphisms it allows and compute $n$ using the formula $n=
2^{(s-r)} \cdot i \cdot h$. Note that $s-r$ is just the number of
irreducible components of $\tilde D$.

\textbf{Example:} We take the diagram of the object $(D;\{A,B\})$
corresponding to a general object of a given stratum, and reduce it
to a diagram of $(\tilde D; \{(\tilde A, H_A), (\tilde B,
H_B)\},H_{AB})$ in the following way: We keep the markings that do
not lie on extremities, and we introduce for every point to which an
extremity is contracted a circle, in the center of which we insert a
dot if the extremity carried two dots, a square if the extremity
carried two squares, and a cross if the extremity carried one square
and one dot. A automorphism must either take all symbols to symbols
of the same kind (i.e. dots to dots, squares to squares, circled
dots to circled dots,...) or it it must take all dots to squares and
vice versa, all circled dots to circled squares and vice versa, and
take circled crosses to circled crosses.

For example, in the case of the stratum $L^+$ we get

\begin{center}
 \begin{tabular}{m{2cm} m{1.5cm} m{2cm}}
  \includegraphics[width=1.6cm]{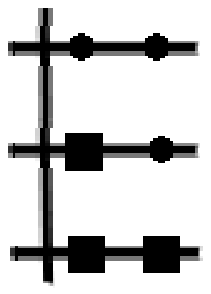} & $\longrightarrow$ & \includegraphics[width=1.6cm]{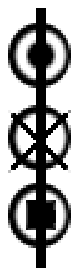}
 \end{tabular}
\end{center}

For $M^+$ we get

\begin{center}
 \begin{tabular}{m{2cm} m{1.5cm} m{2cm}}
  \includegraphics[width=1.6cm]{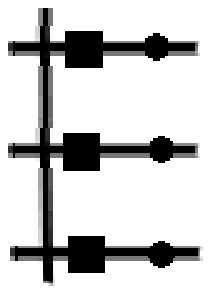} & $\longrightarrow$ & \includegraphics[width=1.6cm]{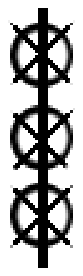}
 \end{tabular}
\end{center}

Now, using Lemma \ref{ap2} (ii), it is clear that $h=2$ for $L^+$
(it is possible to swap the square and the dot), and $h=6$ for
$M^+$. Since both diagrams have only one irreducible component,
$s-r=1$ in both cases. The nonexceptional subcurve of a general
object of $L^+$ has one connected components, so here $i=2^{1-1}=1$,
while for $M^+$ the nonexceptional subcurve has two connected
components, so there $i=2^{2-1}=2$. Putting all this together we get
that the automorphism number $n$ is $4$ for $L^+$ and $24$ for
$M^+$.

The following table contains for each closed stratum of $\bM_2$ the
automorphism numbers of objects corresponding to general points of
the different closed strata of $\bM_2$, and the diagram of the
object corresponding to the preimage of such a general point under
the isomorphism $b: \bM_{0,[6]} \rightarrow \bM_2$, which can be determined easily using the explicit description of 
$b$ in \cite{MR1981190}.  

\begin{center}
\begin{tabular}{|m{1.5cm}|m{2.2cm}|m{2.2cm}|m{3cm}|}
\hline Codim. &\textbf{Stratum} & \textbf{Diagram} & \textbf{Auto. number}\\
\hline
\hline 0 & $S_2$ & \includegraphics[width=1.6cm]{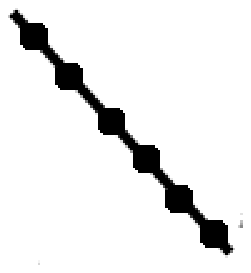} & $2$\\
\hline 1 & $\Delta_0$ & \includegraphics[width=1.6cm]{del_0.eps} & $2$\\
\hline 1 & $\Delta_1$ & \includegraphics[width=1.6cm]{del_1.eps} & $4$\\
\hline 2 & $\Delta_{00}$ & \includegraphics[width=1.6cm]{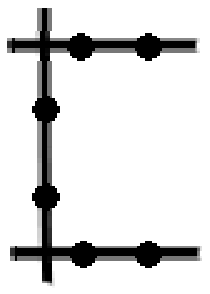}  & $4$\\
\hline 2 & $\Delta_{01}$ & \includegraphics[width=1.6cm]{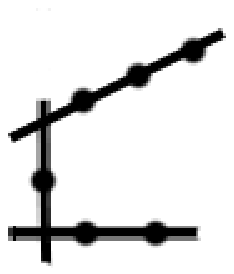} & $4$\\
\hline 2 & $C_{000}$ & \includegraphics[width=1.6cm]{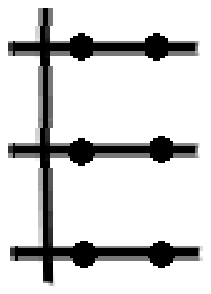} & $12$\\
\hline 2 & $C_{001}$ & \includegraphics[width=1.6cm]{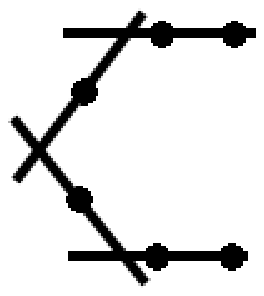} & $8$\\
\hline
\end{tabular}
\end{center}

Next, the same for $S_2^+$, but with an additional column, showing
to which class in $A^*_{\mathbb{Q}} (\bM_2)$ the $Q$-class of each
closed stratum is pushed forward by $\pi_+$. Concerning the
diagram: Here we of course list the diagram belonging of the
preimage of a general point under $a_+: \bM_{0,[3,3]} \rightarrow
\bS_2^+$. Which kind of diagram coresponds to the general spin curve of a stratum can be determined by using that we know such a corespondence already for $\bM_2$, that we know the corespondence for all codimension $1$ Strata from Section \ref{sec1}, and by considdering how their general spin curves can degenerate, and, when in doubt, by counting the degree of the strata over $\bM_2$ and $\bM_{0,[6]}$, like in section \ref{sec1}.
  
\begin{center}
\begin{tabular}{|m{1.5cm}|m{2cm}|m{2cm}|m{2cm}|m{2cm}|}
\hline Codim. &\textbf{Stratum} & \textbf{Diagram} & \textbf{Auto's}& $(\pi_+)_*([...]_Q)$\\
\hline \hline 0 & $\bS_2^+$ & \includegraphics[width=1.6cm]{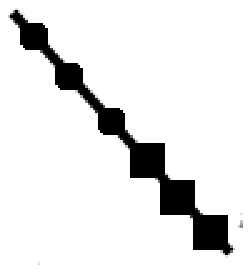} & $2$ & $10 [\bM_2]_Q$ \\
\hline \hline 1& $\Aop$ & \includegraphics[width=1.6cm]{aop.eps} & $2$ & $4 \delta_0$\\
\hline 1& $\Bop$ & \includegraphics[width=1.6cm]{bop.eps} & $2$ & $3 \delta_0$ \\
\hline 1& $\Aip$ & \includegraphics[width=1.6cm]{aip.eps} & $8$ & $\frac{9}{2} \delta_1$\\
\hline 1& $\Bip$ & \includegraphics[width=1.6cm]{bip.eps} & $8$ & $\frac{1}{2} \delta_1$\\
\hline
\end{tabular}

\begin{tabular}{|m{1.5cm}|m{2cm}|m{2cm}|m{2cm}|m{2cm}|}
\hline 2& $C^+$ & \includegraphics[width=1.6cm]{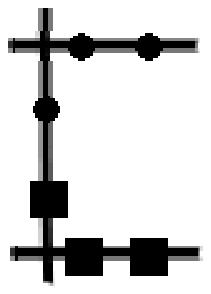} & $4$ & $2 [\Delta_{00}]_Q$\\
\hline 2& $D^+$ & \includegraphics[width=1.6cm]{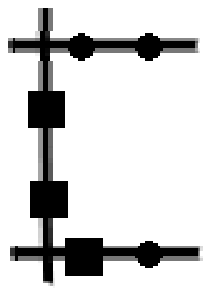} & $2$ & $2 [\Delta_{00}]_Q$\\
\hline 2& $E$ & \includegraphics[width=1.6cm]{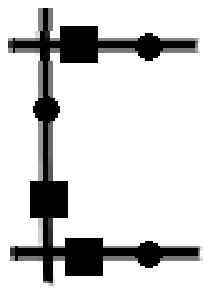} & $4$ & $[\Delta_{00}]_Q$\\
\hline
\end{tabular}

\begin{tabular}{|m{1.5cm}|m{2cm}|m{2cm}|m{2cm}|m{2cm}|}
\hline 2& $X^+$ & \includegraphics[width=1.6cm]{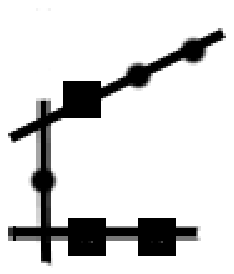} & $8$ & $\frac{3}{2} [\Delta_{01}]_Q$\\
\hline 2& $Y^+$ & \includegraphics[width=1.6cm]{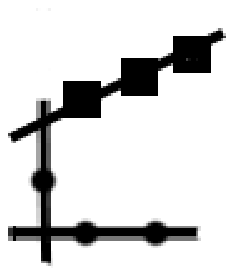} & $8$ & $\frac{1}{2} [\Delta_{01}]_Q$\\
\hline 2& $Z^+$ & \includegraphics[width=1.6cm]{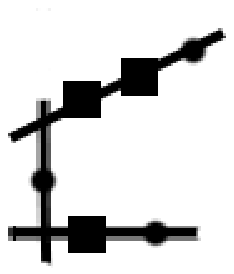} & $8$ & $\frac{3}{2} [\Delta_{01}]_Q$\\
\hline
\end{tabular}

\begin{tabular}{|m{1.5cm}|m{2cm}|m{2cm}|m{2cm}|m{2cm}|}
\hline 3& $L^+$ & \includegraphics[width=1.6cm]{L+.eps} & $4$ & $3 [\Delta_{000}]_Q$ \\
\hline 3& $M$ & \includegraphics[width=1.6cm]{M+.eps} & $24$ & $\frac{1}{2} [\Delta_{000}]_Q$\\
\hline
\end{tabular}

\begin{tabular}{|m{1.5cm}|m{2cm}|m{2cm}|m{2cm}|m{2cm}|}
\hline 3& $Q^+$ & \includegraphics[width=1.6cm]{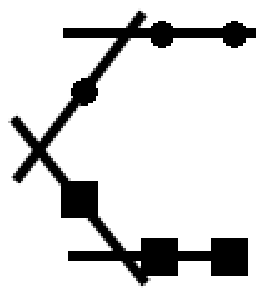} & $16$ & $\frac{1}{2} [\Delta_{001}]_Q$\\
\hline 3& $P^+$ & \includegraphics[width=1.6cm]{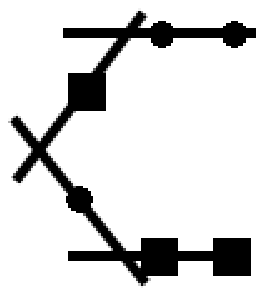} & $16$ & $\frac{1}{2} [\Delta_{001}]_Q$\\
\hline 3& $U^+$ & \includegraphics[width=1.6cm]{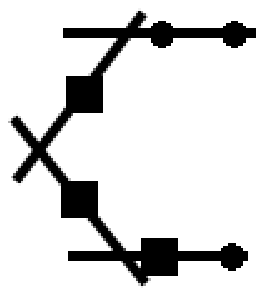} & $8$ & $[\Delta_{001}]_Q$\\
\hline 3& $R$ & \includegraphics[width=1.6cm]{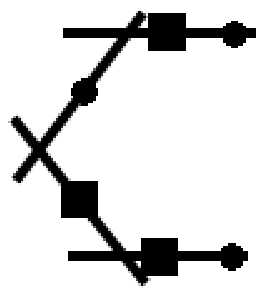} & $16$ & $\frac{1}{2} [\Delta_{001}]_Q$ \\
\hline
\end{tabular}
\end{center}

Next, the same for $\bS_2^-$:

\begin{center}
\begin{tabular}{|m{1.5cm}|m{2cm}|m{2cm}|m{2cm}|m{2cm}|}
\hline Codim. &\textbf{Stratum} & \textbf{Diagram} & \textbf{Auto's}& $(\pi_-)_*([...]_Q)$\\
\hline \hline 0 & $\bS_2^-$ & \includegraphics[width=1.6cm]{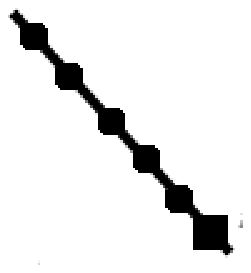} & $2$ & $6 [\bM_2]_Q$ \\
\hline \hline 1& $\Aom$ & \includegraphics[width=1.6cm]{aom.eps} & $2$ & $4 \delta_0$ \\
\hline 1& $\Bom$ & \includegraphics[width=1.6cm]{bom.eps} & $2$ & $\delta_0$\\
\hline 1& $\Aim$ & \includegraphics[width=1.6cm]{aim.eps} & $8$ & $3 \delta_1$\\
\hline
\end{tabular}

\begin{tabular}{|m{1.5cm}|m{2cm}|m{2cm}|m{2cm}|m{2cm}|}
\hline 2& $C^-$ & \includegraphics[width=1.6cm]{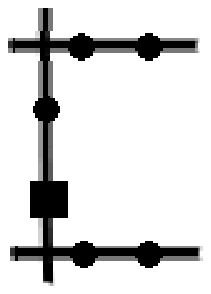} & $2$ & $2 [\Delta_{00}]_Q$ \\
\hline 2& $D^-$ & \includegraphics[width=1.6cm]{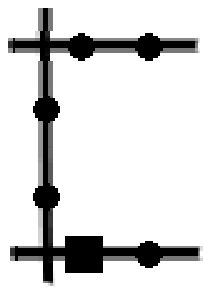} & $2$ & $2 [\Delta_{00}]_Q$\\
\hline
\end{tabular}

\begin{tabular}{|m{1.5cm}|m{2cm}|m{2cm}|m{2cm}|m{2cm}|}
\hline 2& $X^-$ & \includegraphics[width=1.6cm]{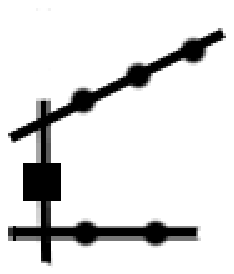} & $8$ & $\frac{1}{2} [\Delta_{01}]_Q$\\
\hline 2& $Y^-$ & \includegraphics[width=1.6cm]{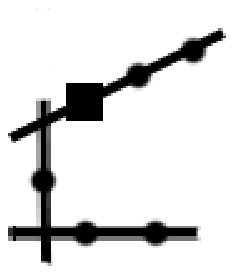} & $8$ & $\frac{3}{2} [\Delta_{01}]_Q$\\
\hline 2& $Z^-$ & \includegraphics[width=1.6cm]{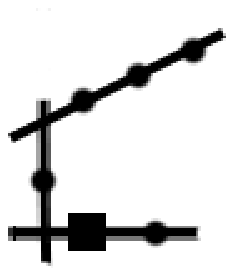} & $8$ & $\frac{1}{2} [\Delta_{01}]_Q$\\
\hline
\end{tabular}

\begin{tabular}{|m{1.5cm}|m{2cm}|m{2cm}|m{2cm}|m{2cm}|}
\hline 3& $L^-$ & \includegraphics[width=1.6cm]{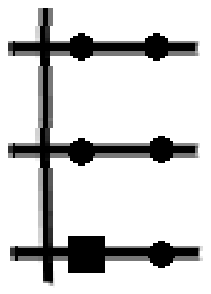} & $4$ & $3 [\Delta_{000}]_Q$\\
\hline 3& $P^-$ & \includegraphics[width=1.6cm]{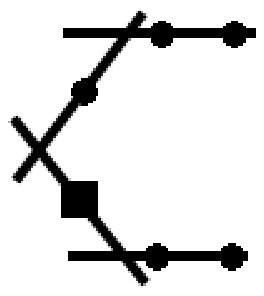} & $8$ & $[\Delta_{001}]_Q$\\
\hline 3& $U^-$ & \includegraphics[width=1.6cm]{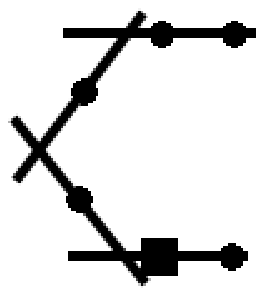} & $8$ & $[\Delta_{001}]_Q$\\
\hline
\end{tabular}
\end{center}

Next, the same for $\bR_2$:

\begin{center}
\begin{tabular}{|m{1.5cm}|m{2cm}|m{2cm}|m{2cm}|m{2cm}|}
\hline Codim. &\textbf{Stratum} & \textbf{Diagram} & \textbf{Auto's}& $(\pi_R)_*([...]_Q)$\\
\hline \hline 0 & $\bR_2$ & \includegraphics[width=1.6cm]{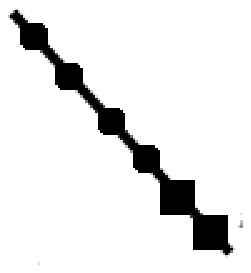} & $2$ & $15 [\bM_2]_Q$ \\
\hline \hline 1 & $\Doi$ & \includegraphics[width=1.6cm]{doi.eps} & $2$ & $6 \delta_0$\\
\hline 1 & $\Doii$ & \includegraphics[width=1.6cm]{doii.eps} & $2$ & $\delta_0$\\
\hline 1 & $\Dor$ & \includegraphics[width=1.6cm]{dor.eps}  & $2$ & $4 \delta_0$ \\
\hline
\end{tabular}

\begin{tabular}{|m{1.5cm}|m{2cm}|m{2cm}|m{2cm}|m{2cm}|}
\hline 1 & $\Di$  & \includegraphics[width=1.6cm]{di.eps} & $4$ & $6 \delta_1$ \\
\hline 1 & $\Dii$ & \includegraphics[width=1.6cm]{dii.eps} & $4$ & $9 \delta_1$\\
\hline
\end{tabular}

\begin{tabular}{|m{1.5cm}|m{2cm}|m{2cm}|m{2cm}|m{2cm}|}
\hline 2 & $E^{\prime,\prime}$ & \includegraphics[width=1.6cm]{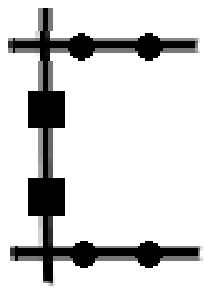} & $4$ & $[\Delta_{00}]_Q$\\
\hline 2& $E^{\prime,\prime \prime}$ & \includegraphics[width=1.6cm]{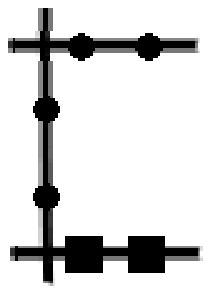} & $2$ & $2 [\Delta_{00}]_Q$\\
\hline 2& $E^{\prime,r}$ & \includegraphics[width=1.6cm]{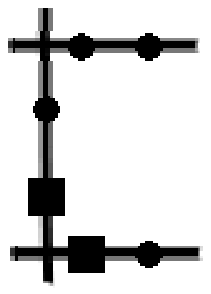} & $2$ & $4 [\Delta_{00}]_Q$\\
\hline 2& $E^{r,r}$ & \includegraphics[width=1.6cm]{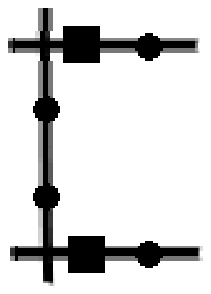} & $4$ & $[\Delta_{00}]_Q$\\
\hline
\end{tabular}

\begin{tabular}{|m{1.5cm}|m{2cm}|m{2cm}|m{2cm}|m{2cm}|}
\hline 2& $F_1'$ & \includegraphics[width=1.6cm]{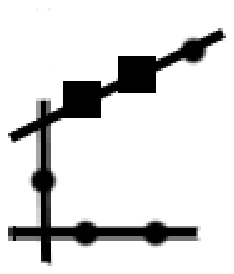} & $4$ & $3 [\Delta_{01}]_Q$\\
\hline 2& $F_1''$ & \includegraphics[width=1.6cm]{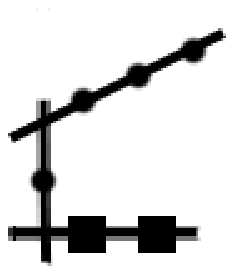} & $4$ & $[\Delta_{01}]_Q$\\
\hline 2& $F_1^r$ & \includegraphics[width=1.6cm]{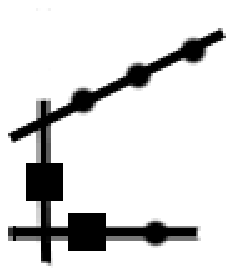} & $4$ & $[\Delta_{01}]_Q$\\
\hline 2& $F_{1:1}'$ & \includegraphics[width=1.6cm]{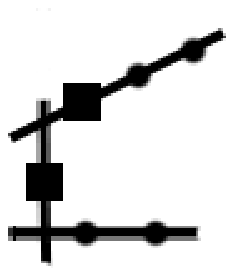} & $4$ & $3 [\Delta_{01}]_Q$\\
\hline 2& $F_{1:1}^r$ & \includegraphics[width=1.6cm]{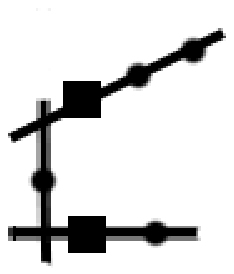} & $4$ & $[\Delta_{01}]_Q$\\
\hline
\end{tabular}

\begin{tabular}{|m{1.5cm}|m{2cm}|m{2cm}|m{2cm}|m{2cm}|}
\hline 3& $G'$ & \includegraphics[width=1.6cm]{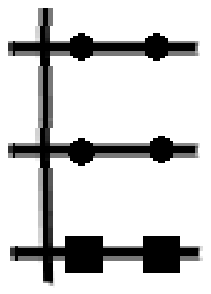} & $4$ & $3 [\Delta_{000}]_Q$\\
\hline 3& $G^r$ & \includegraphics[width=1.6cm]{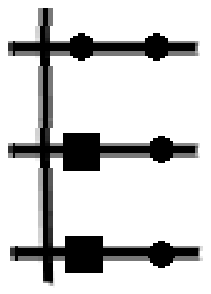} & $4$ & $3 [\Delta_{000}]_Q$\\
\hline
\end{tabular}

\begin{tabular}{|m{1.5cm}|m{2cm}|m{2cm}|m{2cm}|m{2cm}|}
\hline 3& $H_1'$ & \includegraphics[width=1.6cm]{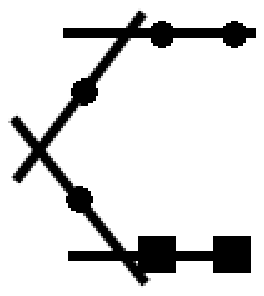} & $4$ & $2 [\Delta_{001}]_Q$\\
\hline 3& $H_1^r$ & \includegraphics[width=1.6cm]{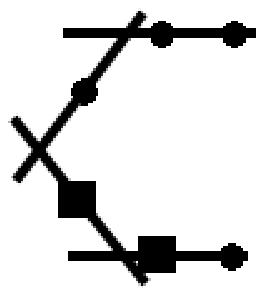} & $4$ & $2 [\Delta_{001}]_Q$\\
\hline 3& $H_{1:1}'$ & \includegraphics[width=1.6cm]{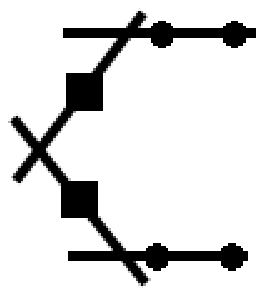} & $8$ & $[\Delta_{001}]_Q$\\
\hline 3& $H_{1:1}^r$ & \includegraphics[width=1.6cm]{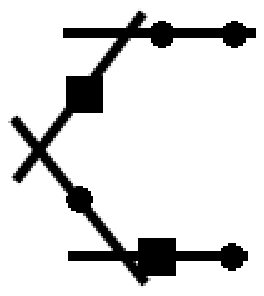} & $4$ & $2 [\Delta_{001}]_Q$\\
\hline 3& $H_{1:1}^{r,r}$ & \includegraphics[width=1.6cm]{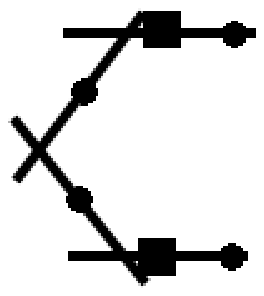} & $8$ & $[\Delta_{001}]_Q$\\
\hline
\end{tabular}
\end{center}



\bibliographystyle{alpha}
\addcontentsline{toc}{section}{References}
\bibliography{BibSeb}

\end{document}